\documentclass[11pt]{article}

\usepackage[utf8]{inputenc}
\usepackage{preamble}

\usepackage{tikz}
\usepackage{circuitikz}
\usetikzlibrary{calc,shapes,arrows,patterns,decorations.pathmorphing,decorations.markings,arrows.meta}

\usetikzlibrary{external}
\title{\textbf{Onset and stabilization of delay-induced instabilities in piezoelectric digital vibration absorbers}}
\author{G. Raze$^{1}$, J. Dietrich$^{1}$ and G. Kerschen$^{1}$\\
        \small $^{1}$Department of Aerospace and Mechanical Engineering, \\
        \small University of Liège, 4000 Liège, Belgium
        }
\date{}

\begin{document}

\maketitle

\abstractn{The stability of a piezoelectric structure controlled by a digital vibration absorber emulating a shunt circuit is investigated in this work. The formalism of feedback control theory is used to demonstrate that systems with a low electromechanical coupling are prone to delay-induced instabilities entailed by the sampling procedure of the digital unit. An explicit relation is derived between the effective electromechanical coupling factor and the maximum sampling period guaranteeing a stable controlled system. Since this sampling period may be impractically small, a simple modification procedure of the emulated admittance of the shunt circuit is proposed in order to counteract the effect of delays by anticipation. The theoretical developments are experimentally validated on a clamped-free piezoelectric beam.}

\keywords{vibration mitigation, piezoelectric shunt damping, digital vibration absorber, delay-induced instabilities, stabilization}

 \section{Introduction}
    
        Engineering structures from various disciplines tend to be lighter or more slender, which usually goes along with smaller structural damping and increased susceptibility to vibrations, threatening structural integrity. Piezoelectric shunt damping is often considered as one potential solution to this issue. It was originally proposed by \cite{Forward1979}, and formalized by \cite{Hagood1991}. The working principle of piezoelectric shunt damping is based on the capability of piezoelectric transducers to convert a part of their mechanical energy into electrical energy. The latter can be dissipated by connecting a shunt circuit to the electrodes of the transducer. A common type of shunt is a resonant one, composed of a resistor and an inductor, arranged either in series or in parallel. The realization of this circuit may be challenging for several reasons. First, the required inductance may be impractically large. Second, the performance of the piezoelectric shunt is highly sensitive to the values of the electrical components. Any misevaluation or time variation of the system characteristics will result in sub-optimal performance, rectified by time-consuming manual modifications of the electrical parameters. 
	
	\cite{Fleming2000,Fleming2004} introduced the concept of synthetic impedance as an alternative solution. The combination of a digital signal processor with a current source makes the realization of an arbitrary impedance possible. The synthetic impedance is an attractive option to realize shunt damping circuits owing to its versatility. This nonetheless comes at the expense of the need for powering the digital unit and its associated electronics. Since it was proposed, the application of piezoelectric shunt damping with a digital vibration absorber (DVA) was used in several works. \cite{Fleming2003a} and \cite{Pliva2007} developed architectures using pulse width modulation in order to simplify the driving electronics. \cite{Niederberger2004} implemented an adaptive impedance with a DVA in order to improve the robustness of the control system.  \cite{Giorgio2009} and \cite{Rosi2010} used digital controllers to validate their theoretical developments on piezoelectric damping with electrical networks. \cite{Matten2014}, \cite{Necasek2016,Necasek2017} and \cite{Silva2018t} investigated various electronic architectures to implement a DVA, and discussed how to set up its analog and digital parts. \cite{DalBo2019} configured a digital unit to realize vibration absorbers with swept and switched characteristics. Recently, this concept was applied to metamaterials by \cite{Sugino2018,Sugino2020} and \cite{Yi2020}, and to nonlinear shunt damping \citep{Raze2020c}.	
		
	In the active control terminology, the DVA is equivalent to a control system with a self-sensing actuator, and it implements a passive control law. From a theoretical standpoint, a passive control law features unconditional stability of the controlled system \citep{Moheimani2003}. If the problem is cast into a feedback control one, the system exhibits an infinite gain margin but a finite phase margin. Because a digital unit needs to sample the signals it is working with, unavoidable delays occur in the control loop. These delays introduce a phase lag which may destabilize the controlled system if they are too large. \cite{Necasek2016} and \cite{Sugino2018} pinpointed the fact that in some cases a DVA needs to have a sampling frequency much higher than the typical frequencies of interest. The authors also noted that delay-induced instabilities may arise when using a DVA for surprisingly small sampling periods, in spite of the passivity of the control law~\citep{Raze2019}.
	
	In this work, novel and ready-to-use formulas are provided to determine whether delay-induced instabilities can be an issue, and how to solve this issue if need be. Specifically, the purpose of this paper is to evidence 1. why a DVA may need a high sampling frequency for stability, 2. how delay-induced instabilities may arise and 3. how to counteract them.  After reviewing the basics of piezoelectric shunt damping with a DVA in Section~\ref{sec:piezoShuntDVA}, the possibility for instabilities of the control system is investigated in Section~\ref{sec:instabilities}. The problem is cast as a feedback control one, and a relation between the effective electromechanical coupling factor and the phase margin is derived. Values of the maximum sampling period under which the system remains stable are then deduced. A procedure to stabilize the controlled system is proposed in Section~\ref{sec:stabilization}. Upon applying this procedure, larger sampling periods may be used for the digital unit without jeopardizing stability, which is generally advantageous. The findings are experimentally validated with a piezoelectric clamped-free beam in Section~\ref{sec:experiments}. The conclusions of this work are finally reported in Section~\ref{sec:conclusion}. In comparison to the conference paper~\citep{Raze2019}, this works performs an in-depth investigation of the controlled system's dynamics, links the delay-induced instabilities to the electromechanical coupling and improves the stabilization procedure proposed therein.

    \section{Piezoelectric shunt damping with a digital vibration absorber}
    \label{sec:piezoShuntDVA}

        A single-degree-of-freedom structure to which a piezoelectric transducer is bonded is considered. The structure is excited by an external force $f$ and responds with a displacement $x$. $V$ and $\dot{q}$ (where an upper dot denotes time derivation) are the voltage across the electrodes of the transducer and the current flowing through them, respectively. The governing equations of the piezoelectric structure read
  		\begin{equation}
			\left\{
			\begin{array}{l}
				m\ddot{x} + k_{oc}x - \theta_p q = f \\
				V = \theta_p x - \dfrac{1}{C_p^\varepsilon}q
			\end{array}
			\right. ,
			\label{eq:EoM-Structure}
		\end{equation}
		where $m$ and $k_{oc}$ are the mass and stiffness of the structure when the transducer is open-circuited, respectively, $\theta_p$ is a piezoelectric coupling coefficient and $C_p^\varepsilon$ is the piezoelectric capacitance at constant strain. The resonance frequency of the structure when the transducer is open-circuited ($q=0$) is given by
		\begin{equation}
		    \omega_{oc} = \sqrt{\dfrac{k_{oc}}{m}},
		    \label{eq:ocFrequency}
		\end{equation}
		and when the transducer is short-circuited ($V=0$), the stiffness of the structure changes to $k_{sc}$, and the short-circuit resonance frequency can be found as
		\begin{equation}
		    \omega_{sc} = \sqrt{\dfrac{k_{sc}}{m}} = \sqrt{\dfrac{k_{oc}-\theta_p^2C_p^\varepsilon}{m}}.
		    \label{eq:scFrequency}
		\end{equation}
		The electromechanical coupling between the transducer and the structure can be assessed from these two frequencies with the effective electromechanical coupling factor (EEMCF) $K_c$, defined by
		\begin{equation}
		    K_c^2 = \dfrac{\omega_{oc}^2-\omega_{sc}^2}{\omega_{sc}^2}.
		    \label{eq:EEMCF}
		\end{equation}
		
		Upon connecting a series RL shunt to the transducer, the equations of the system become
		\begin{equation}
			\left\{
			\begin{array}{l}
				m\ddot{x} + k_{oc}x - \theta_p q = f \\
				L\ddot{q} + R\dot{q} + \dfrac{1}{C_p^\varepsilon}q - \theta_p x = 0
			\end{array}
			\right. ,
			\label{eq:EoM-SeriesRL}
		\end{equation}
		where $L$ and $R$ are the inductance and resistance of the shunt, respectively. They can be optimized to minimize the maximum vibratory amplitude of the structure  \citep{Soltani2014a,Ikegame2019}. Introducing an intermediate parameter
		\begin{equation}
		    r = \frac{\sqrt{64-16 K_c^2-26 K_c^4}-K_c^2}{8},
		    \label{eq:Seriesr}
		\end{equation}
		The optimal inductance and resistance are
		\begin{equation}
		    L = \frac{4 K_c^2+4}{3 K_c^2-4 r+8}\frac{1}{\omega_{oc}^2C_p^\varepsilon} = \dfrac{1}{\delta^2(K_c)\omega_{oc}^2C_p^\varepsilon}
		    \label{eq:SeriesLopt}
		\end{equation}
		and
		\begin{equation}
		    R = \frac{2 \sqrt{2\left(K_c^2+1\right) \left[27 K_c^4+K_c^2 (80-48 r)-64 (r-1)\right]}}{\left(5 K_c^2+8\right) \sqrt{3K_c^2-4 r+8}}\frac{1}{\omega_{oc}C_p^\varepsilon}=\frac{2\zeta(K_c)}{\delta(K_c)\omega_{oc}C_p^\varepsilon},
		    \label{eq:SeriesRopt}
		\end{equation}
		respectively, where $\delta$ and $\zeta$ represent electrical frequency and damping ratios, respectively. With this tuning, the frequency response function (FRF) of the controlled structure exhibits two peaks of identical amplitude, generally much lower than the maximum amplitude of the uncontrolled structure.
		
\begin{figure}[!ht]
        		\centering
        		\scalebox{0.7}{
        		\includegraphics[scale=1]{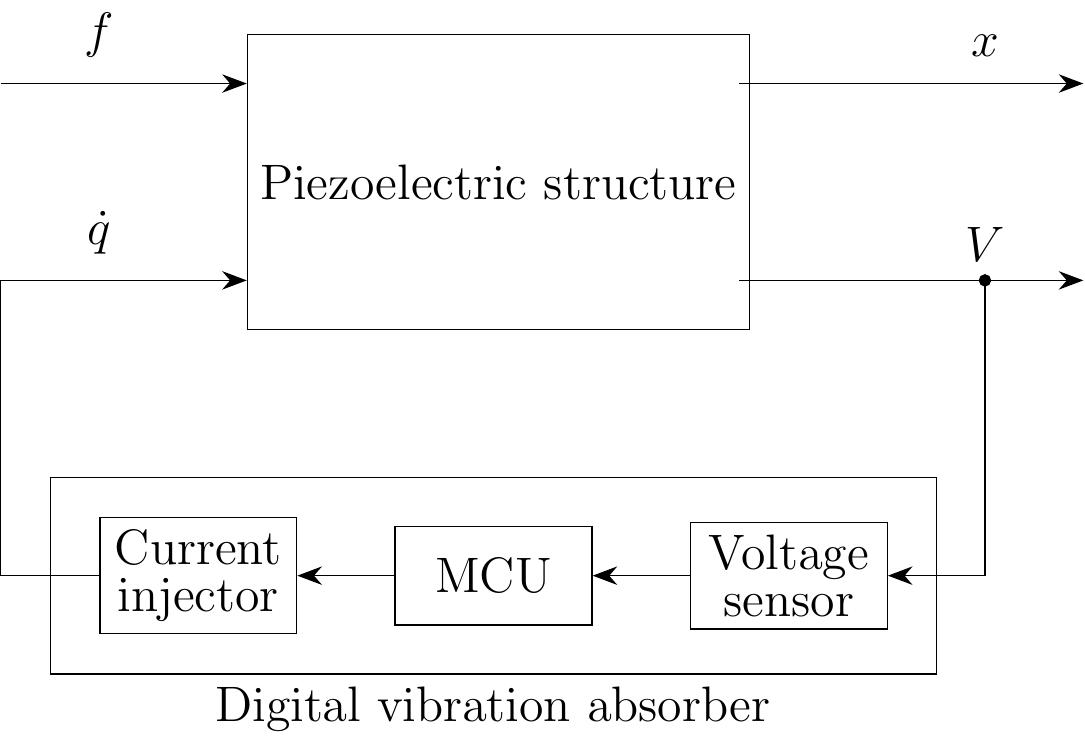}
        		}
        		\caption{General working principle of the DVA.}
        		\label{fig:DAPrinciples}
        	\end{figure}

		With typical values of EEMCF and piezoelectric capacitance, the required inductance is usually large, which makes the practical realization of a RL shunt difficult. This is one of the reasons that motivated the use of DVAs. This concept, proposed by \cite{Fleming2000}, consists in replacing an analog shunt by a digital electronic circuit mimicking the former. It does so by measuring the voltage of a piezoelectric transducer and by injecting the desired current back into it. The voltage-to-current transfer function, i.e. the admittance, is programmed into a microcontroller unit (MCU). The general working principle of this approach is schematized in Fig.~\ref{fig:DAPrinciples}.

    \section{Delay-induced instabilities}
    \label{sec:instabilities}
        
        The purpose of this section is to study the stability of a system composed of a piezoelectric structure controlled by a DVA. Using the theory of feedback control \citep{Franklin2015}, it is demonstrated that the system can be sensitive to delays induced by the sampling made by the digital unit.
    
        \subsection{Delays induced by the sampling procedure}
         \label{ssec:samplingDelays}
         
\begin{figure*}[!ht]
    			\centering
    			\scalebox{0.8}{
        		\includegraphics[scale=1]{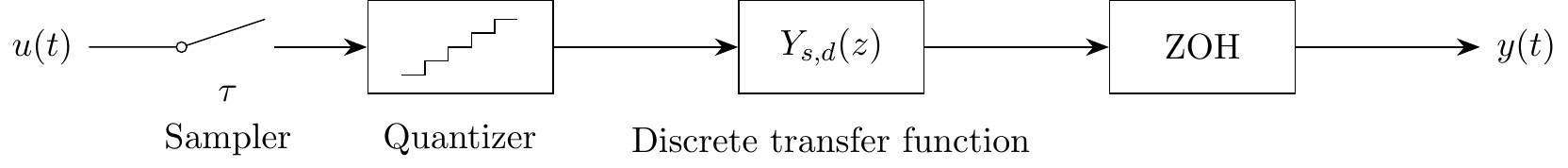}
    			}
    			\caption{Block diagram representation of the input/output relation in a digital system.}
    			\label{fig:mcuSchematics}
    		\end{figure*}

            Fig.~\ref{fig:mcuSchematics} depicts a schematic representation of the process undergone by an input signal $u(t)$ (typically, the piezoelectric voltage) to be transformed into an output signal $y(t)$ (typically, the piezoelectric current) by a digital unit \citep{Franklin1998}.	
    	
            A sample-and-hold circuit holds the input signal $u(t)$ constant at specific times, multiples of the sampling period $\tau$. The signal is then quantized, and the MCU operates on it to emulate the desired admittance. This signal being discrete, a discrete input-output transfer function must thus be employed. Tustin's method \citep{Tustin1947} is used to discretize the continuous transfer function to be emulated. If the latter is given by $Y_s(s)$, a discrete $z$-transform $Y_{s,d}(z)$ is obtained by substituting the $s$ variable by a bilinear function of the complex variable $z$ as
        	\begin{equation}
        		Y_{s,d}(z) = \left. Y_s(s) \right|_{s=\frac{2}{\tau}\frac{z-1}{z+1}}.
        	\end{equation} 
        	This transfer function can then readily be translated into a recurrence equation \citep{Franklin1998}. The resulting output signal is also a discrete signal. It is applied to the continuous system by holding its value constant for the sample interval by a zero-order hold (ZOH), which keeps the output signal constant over a whole sampling interval.
        	
        	If the MCU operates at a high enough clock frequency relative to the sampling frequency, it may be considered that the digitization of the input signal and computation of the output signal occur instantaneously at each sampling time, i.e., latency is neglected. The differences between the continuous transfer function and the digital one then principally comes from the delay brought by the ZOH, as well as the frequency warping stemming from the discretization of the transfer function.
        	
\begin{figure}[!ht]
                \centering
                \includegraphics[width=0.45\textwidth]{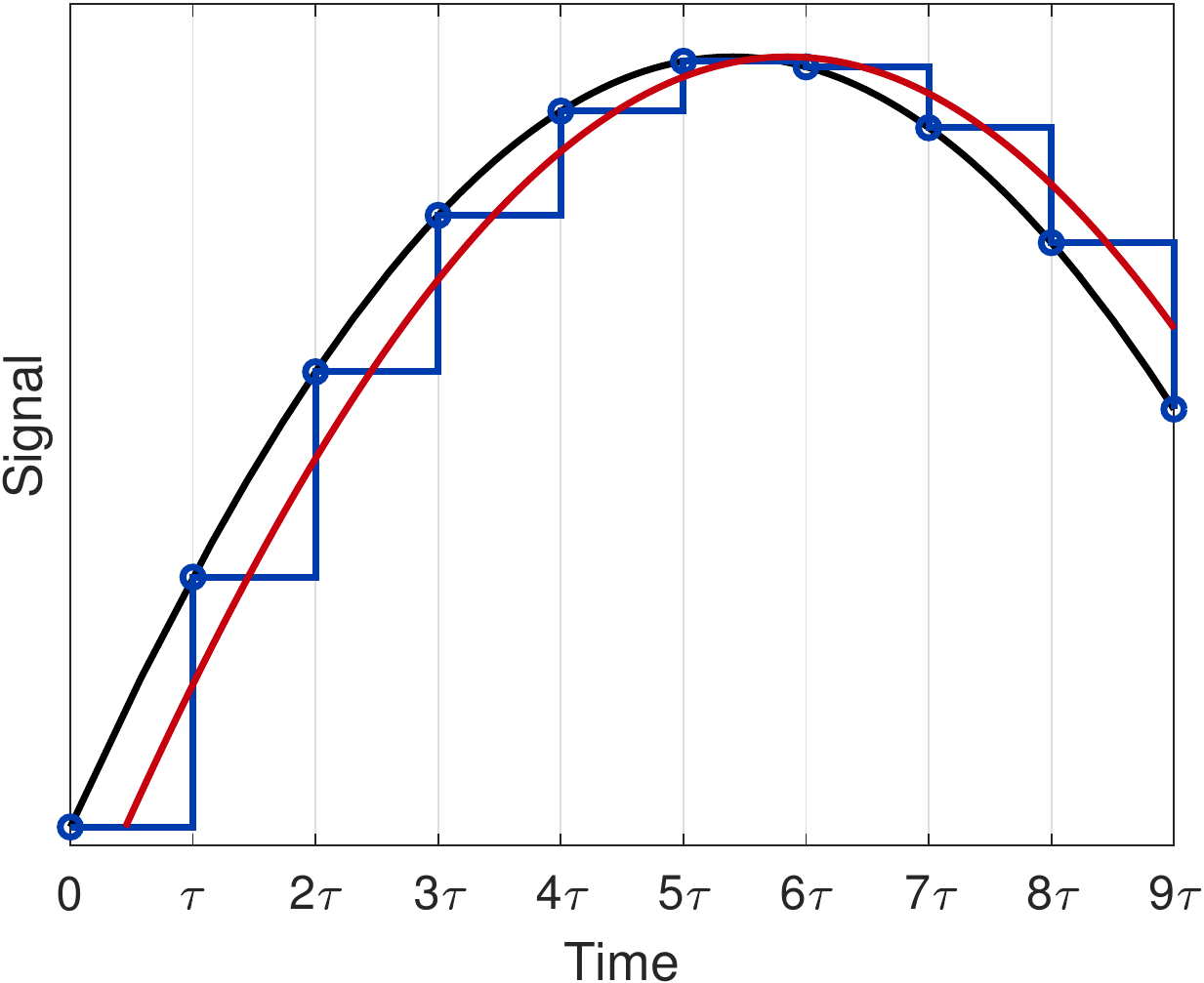}
                \caption{Signals in the MCU when it operates as a simple unity gain: input signal (\rule[0.2em]{1em}{0.2em}), output signal (\textcolor{col1}{\rule[0.2em]{1em}{0.2em}}) and continuous average of the output signal (\textcolor{col2}{\rule[0.2em]{1em}{0.2em}}).}
                \label{fig:zoh}
            \end{figure}

        	An example relating the input and output signals when the MCU implements a simple unity gain ($Y_s=1$) is featured in Fig.~\ref{fig:zoh}. In this case, a continuous average of this output signal looks identical to the input signal but delayed by $\tau/2$.
    
    	\subsection{Open-loop analysis}
			
			%After representing the action of the shunt circuit as a feedback, it is shown that the phase margin can be small for systems with small electromechanical coupling, which makes them prone to delay-induced instabilities because delays induce additional phase lags. Before moving on to a closed-loop analysis, a physical explanation on why the system is sensitive to phase lags is provided, showing that delays make the digital admittance non-passive.
			
			\subsubsection{Open-loop transfer function}
						
\begin{figure}[!ht]
					\centering
					\begin{subfigure}{.45\textwidth}				
						\centering
					\scalebox{0.8}{
						\includegraphics[scale=1]{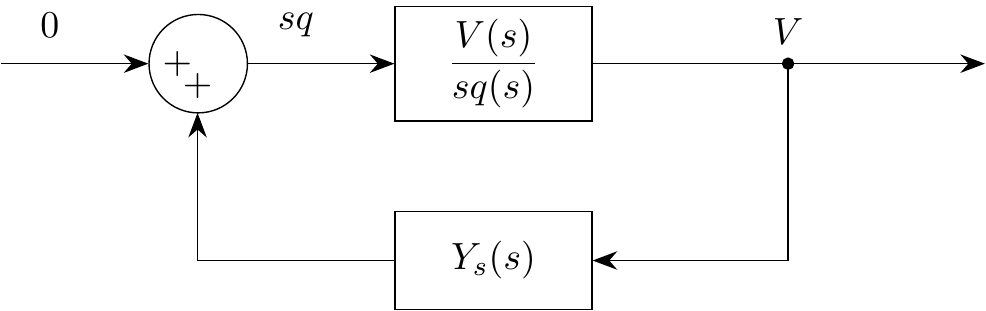}
					}
					\caption{}
					\label{fig:feedbackShunt}
					\end{subfigure}

					\begin{subfigure}{.45\textwidth}	
					\centering
					\scalebox{0.8}{
					
        		\includegraphics[scale=1]{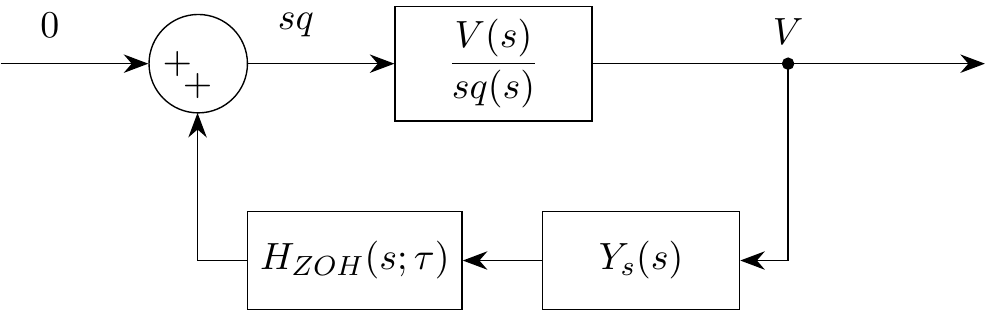}
					}
					\caption{}
					\label{fig:feedbackDelayedShunt}
					\end{subfigure}
					\caption{Block diagram representation of the nominal~\subref{fig:feedbackShunt} and delayed~\subref{fig:feedbackDelayedShunt} controlled systems.}
				\end{figure}

				The stability of the nominal controlled system (i.e., without delays) is conditioned upon that of the unforced system ($f=0$). In this case, the Laplace transform of Eq.~\eqref{eq:EoM-Structure} gives a relation between the piezoelectric voltage and charge (using Eqs.~\eqref{eq:ocFrequency} and \eqref{eq:scFrequency})
				\begin{equation}
				    \dfrac{q}{V} = -C_p^\varepsilon \left(1-\dfrac{C_p^\varepsilon\theta_p^2}{ms^2+k_{oc}}\right)^{-1} = -C_p^\varepsilon \left(1-\dfrac{\omega_{oc}^2-\omega_{sc}^2}{s^2+\omega_{oc}^2}\right)^{-1}  = -C_p^\varepsilon \dfrac{s^2+\omega_{oc}^2}{s^2+\omega_{sc}^2}.
				    \label{eq:dynamicCapacitance}
				\end{equation}
				 This transfer function is usually called the dynamic capacitance \citep{Preumont2011}. Moreover, connecting an admittance $Y_s(s)$ to the electrodes of the transducer imposes the relation
				\begin{equation}
					sq = Y_s(s) V = \dfrac{1}{Ls+R} V
					\label{eq:shuntAdmittanceRelation}
				\end{equation}
				in the case of a series RL shunt. This suggests that the dynamics of the unforced controlled system may be represented with the feedback diagram depicted in Fig.~\ref{fig:feedbackShunt}.
				
				Using Eqs.~\eqref{eq:dynamicCapacitance}  and \eqref{eq:shuntAdmittanceRelation}, one may form the open-loop transfer function $H$ associated with this feedback loop
				\begin{equation}
					H(s) = -\dfrac{V(s)}{sq(s)}Y_s(s)= \dfrac{1}{C_p^\varepsilon} \dfrac{s^2 + \omega_{sc}^2}{s^2 + \omega_{oc}^2} \dfrac{1}{Ls^2+Rs}
					\label{eq:characteristicOL}
				\end{equation}
				such that the poles of the closed-loop system may be found by solving the characteristic equation
				\begin{equation}
					1 + H(s) = 0.
					\label{eq:characteristicOL2}
				\end{equation}

				By normalizing the Laplace variable with the short-circuit resonance frequency
				\begin{equation}
					\overline{s} = \dfrac{s}{\omega_{sc}}
				\end{equation}
				and using Eqs.~\eqref{eq:EEMCF},  \eqref{eq:SeriesLopt}, \eqref{eq:SeriesRopt} and \eqref{eq:characteristicOL},
				\begin{equation}
					H(s) = H\left(\omega_{sc}\overline{s}\right) = \dfrac{\overline{s}^2 + 1}{\overline{s}^2 + 1+K_c^2}\dfrac{1}{\dfrac{1}{\left(1+K_c^2\right)\delta^2(K_c)}\overline{s}^2 + \dfrac{2\zeta(K_c)}{\delta(K_c)\sqrt{1+K_c^2}}\overline{s}},
					\label{eq:characteristicNormalizedOL}
				\end{equation}
				the coefficients of the open-loop transfer function depend only on the EEMCF (since $\delta$ and $\zeta$ are sole functions of it). This parameter is thus expected to play an important role in stability.
					
\begin{figure}[!ht]
					\begin{subfigure}{.45\textwidth}				
						\centering
						\includegraphics[width=\textwidth]{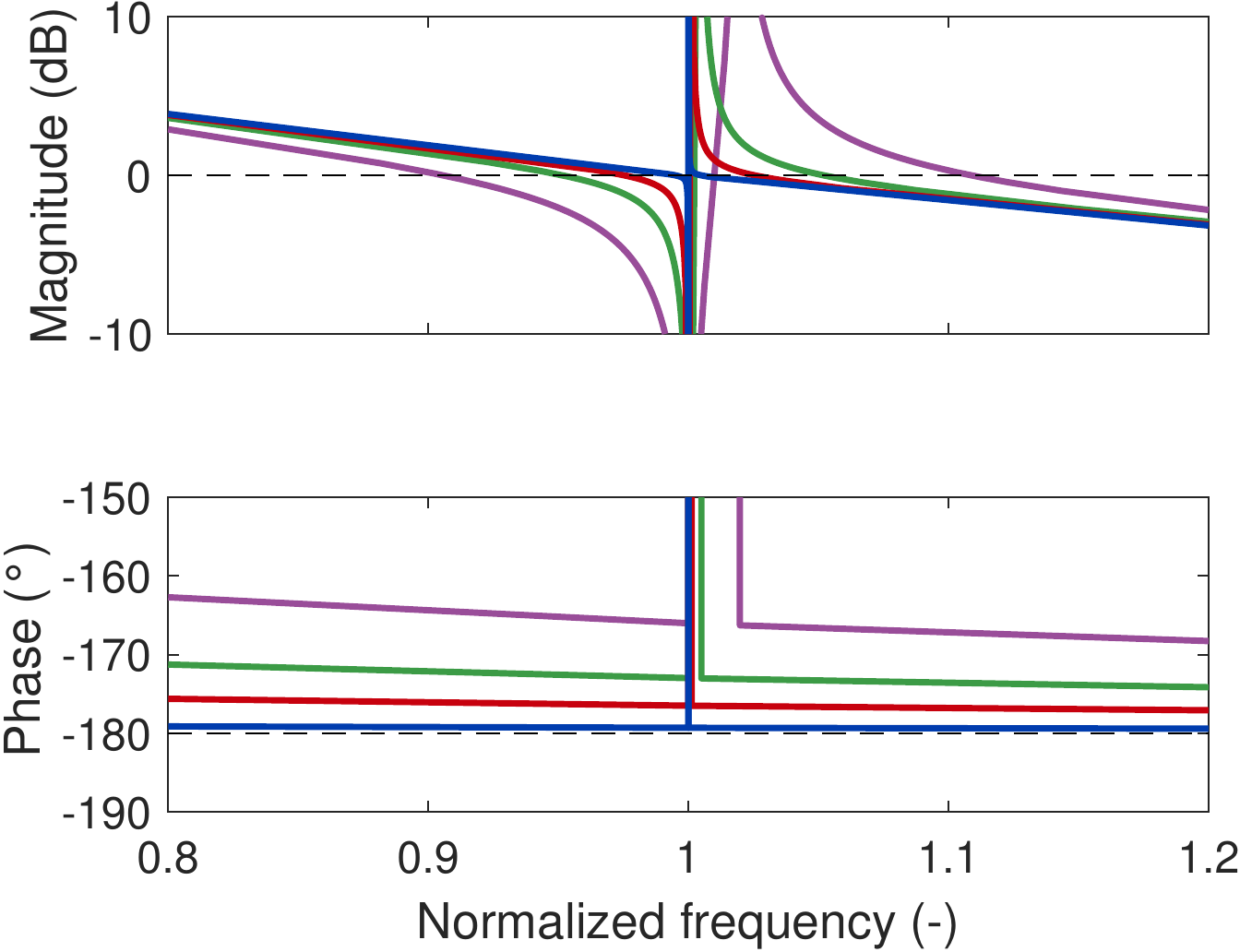}
						\caption{}
		  					\label{sfig:openLoopTF}
					\end{subfigure}
					\hspace{0.05\textwidth}
					\begin{subfigure}{.45\textwidth}				
						\centering
						\includegraphics[width=\textwidth]{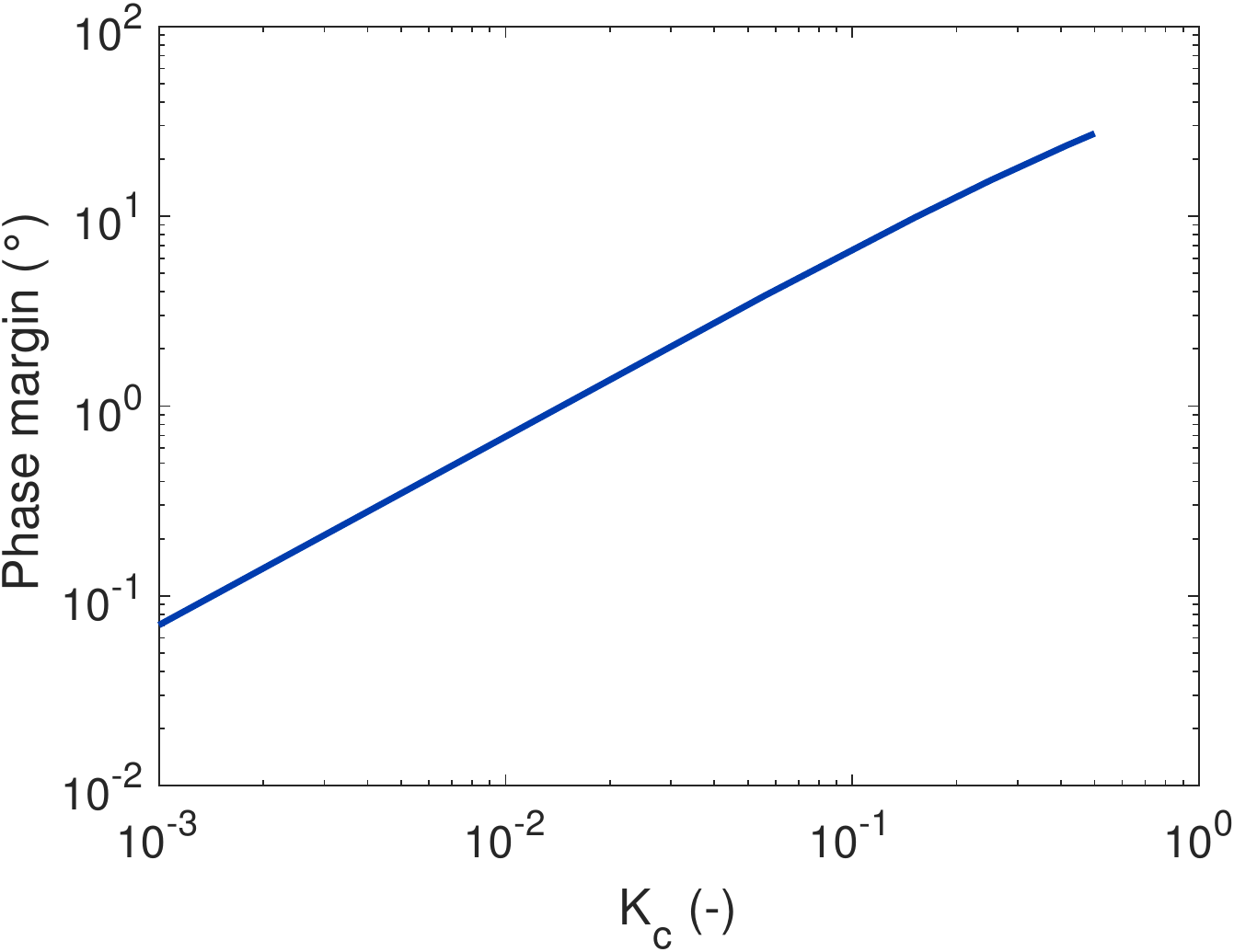}
						\caption{}
		  					\label{sfig:phaseMargin}
					\end{subfigure}
					\caption{Bode plot of the open-loop transfer function~\subref{sfig:openLoopTF}: $K_c=0.01$ (\textcolor{col1}{\rule[0.2em]{1em}{0.2em}}), $K_c=0.05$ (\textcolor{col2}{\rule[0.2em]{1em}{0.2em}}), $K_c=0.1$ (\textcolor{col3}{\rule[0.2em]{1em}{0.2em}}) and $K_c=0.2$ (\textcolor{col4}{\rule[0.2em]{1em}{0.2em}}) ; phase margin as a function of $K_c$~\subref{sfig:phaseMargin}.}
					\label{fig:openLoopTF}
				\end{figure}

				Fig.~\ref{fig:openLoopTF}\subref{sfig:openLoopTF} features Bode plots of the open-loop transfer function given in Eq.~\eqref{eq:characteristicNormalizedOL} for various values of $K_c$ around the short- and open-circuit resonance frequencies. The system has an infinite gain margin because the phase never crosses -180$^\circ$. There are three gain crossover frequencies, and the phase margin is calculated at the highest one (which also corresponds to the lowest phase margin). The phase margin decreases with $K_c$. This trend is also highlighted in Fig.~\ref{fig:openLoopTF}\subref{sfig:phaseMargin}. %\reviewAddition{The phase margin can be small for systems with small electromechanical coupling, which makes them prone to delay-induced instabilities because delays induce additional phase lags.}
				%The EEMCF of a piezoelectric structure is typically a small number, and thus the system may have a small phase margin. This in turn makes it sensitive to time delays in the loop.
				
				\subsubsection{Destabilization mechanism}
				
				In order to intuitively understand why delays can cause instabilities, a simple model is now introduced. As seen in Section~\ref{ssec:samplingDelays}, the delays imparted by the sampling procedure can be modeled as a pure time delay $\tau/2$ \citep{Franklin1998}. In the case of a series RL shunt, the piezoelectric charge and voltage are thus linked by
				\begin{equation}
					L\ddot{q}(t) + R\dot{q}(t) = V\left(t-\dfrac{\tau}{2}\right).
				\end{equation}
				Taking the Laplace transform of this equation yields
				\begin{equation}
					q = \dfrac{e^{-\frac{s\tau}{2}}}{Ls^2 + Rs} V = \dfrac{e^{-\frac{s\tau}{2}}}{s}Y_s(s) = \dfrac{1}{s}Y_d(s),
					\label{eq:simplifiedDelayedAdmittance}
				\end{equation}
				where $Y_s$ is the nominal shunt admittance and $Y_d$ is an equivalent delayed admittance. In order to see how these two quantities differ with a simple exposition, the formulas from \cite{Thomas2012} (which are a linearization of Eqs.~\eqref{eq:SeriesLopt} and \eqref{eq:SeriesRopt} with respect to $K_c$) are used to tune the inductance and the resistance.
				\begin{equation}
					L = \dfrac{1}{C_p^\varepsilon \omega_{oc}^2}, \qquad \qquad R = \sqrt{\dfrac{3}{2}}\dfrac{K_c}{\omega_{oc}C_p^\varepsilon}.
				\end{equation}
				The admittance of the shunt evaluated at $\omega_{oc}$ is thus
				\begin{equation}
					Y_{s}(j\omega_{oc}) = \dfrac{1}{j\omega_{oc}L + R} = \dfrac{C_p^\varepsilon \omega_{oc}}{j + \sqrt{\dfrac{3}{2}}K_c}  = \dfrac{C_p^\varepsilon \omega_{oc}}{1 + \dfrac{3}{2}K_c^2}\left(\sqrt{\dfrac{3}{2}}K_c - j\right).
					\label{eq:nominalAdmittancewoc}
				\end{equation}				
				An important feature of this admittance is that it has a positive real part. Indeed, the average power dissipated across an admittance $Y$ is
				\begin{equation}
					P = \dfrac{1}{2}\Re \left\{ V^* I\right\}= \dfrac{1}{2}\Re \left\{ V^* Y V\right\} =\dfrac{1}{2}\Re \left\{ Y \right\}\left| V \right|^2  
				\end{equation}
				and must be positive for a passive circuit, because it dissipates true power ($\Re$ denotes the operator giving the real part of a complex number and superscript $*$ denotes complex conjugate). Another important feature is that since $K_c \ll 1$, this real part is much lower than the absolute value of the imaginary part.
				
\begin{figure}[!ht]
				\centering
				\scalebox{1.}{
        		\includegraphics[scale=1]{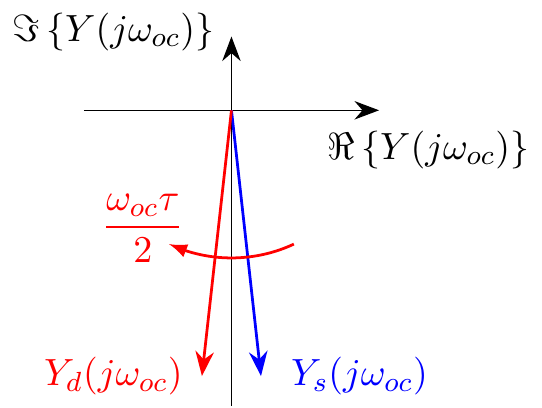}
						}  
					\caption{Representation of the admittance in the complex plane.}
					\label{fig:delayedAdmittance}
				\end{figure}

				The nominal admittance $Y_s(j\omega_{oc})$ given by Eq.~\eqref{eq:nominalAdmittancewoc} is plotted in Fig.~\ref{fig:delayedAdmittance} (where $\Im$ is the operator giving the imaginary part of a complex number). Using Eq.~\eqref{eq:simplifiedDelayedAdmittance}, the delayed admittance $Y_d(j\omega_{oc})$ can be obtained through a clockwise rotation of angle $\omega_{oc} \tau/2$ of the complex vector $Y_s(j\omega_{oc})$. If this angle (i.e., the delay $\tau$) is large enough, this rotation may result in a delayed admittance with a negative real part. This entails the generation of true power within the delayed admittance which may be transmitted to the structure and potentially cause the instability of the closed-loop system. Because the nominal admittance is almost aligned with the imaginary axis for systems with small EEMCFs, this can happen for relatively small delays. 
				
				More quantitatively, using Eqs.~\eqref{eq:simplifiedDelayedAdmittance} and \eqref{eq:nominalAdmittancewoc}, the delayed admittance is given by
				\begin{equation}
					Y_d(j\omega_{oc})  = \dfrac{C_p^\varepsilon \omega_{oc}}{1+\dfrac{3}{2}K_c^2} \left[ \sqrt{\dfrac{3}{2}}K_c \cos \left(\dfrac{\omega_{oc}\tau}{2}\right) - \sin \left(\dfrac{\omega_{oc}\tau}{2}\right) \right.  \left.- j \left(\cos \left(\dfrac{\omega_{oc}\tau}{2}\right) + \sqrt{\dfrac{3}{2}}K_c  \sin \left(\dfrac{\omega_{oc}\tau}{2}\right) \right)\right],
				\end{equation}						
				whose real part becomes negative when
				\begin{equation}
					\tau = \dfrac{2}{\omega_{oc}}\arctan \left( \sqrt{\dfrac{3}{2}} K_c \right) = \dfrac{1}{\omega_{oc}} \left(\sqrt{6} K_c + O(K_c^3)\right) = \dfrac{1}{\omega_{sc}} \left(\sqrt{6} K_c + O(K_c^3)\right) .
					\label{eq:tauNonPassive}
				\end{equation}			
				
		\subsection{Closed-loop analysis}
		
			%The poles of the closed-loop system can be computed with a model accounting for the delays of the ZOH. Root loci given by the variation of the time delay $\tau$ are traced, and the destabilizing effect of the delays can be observed. In order to quantify their importance, approximations of the model are made in order to determine the value of the sampling period $\tau_c$ above which the closed-loop system becomes unstable.
		
			\subsubsection{Characteristic equation}	
				
				The closed-loop system when the structure is controlled by a DVA is represented in Fig.~\ref{fig:feedbackDelayedShunt}. Delays are introduced in the system by the ZOH. 				
				
				%Sampling makes the system time variant. However, if the signals are band-limited (i.e., their frequency content beyond the frequency $\pi / \tau$ is zero so as to respect the conditions for the Nyquist-Shannon sampling theorem), this time-varying character may be neglected. 
				Assuming the output of the ZOH is dominated by the fundamental harmonic of the frequency it is subject to, an equivalent continuous transfer function can be shown to be \citep{Franklin1998}
				\begin{equation}
					H_{ZOH}(s;\tau) = \dfrac{1-e^{-\tau s}}{\tau s}.	
					\label{eq:zoh}
				\end{equation}
				Based on Fig.~\ref{fig:feedbackDelayedShunt}, the characteristic equation is then
				\begin{equation}
					1 - \dfrac{V(s)}{s q(s)} Y_s(s) H_{ZOH}(s;\tau) = 1 + H(s)\dfrac{1-e^{-\tau s}}{\tau s} = 0.
					\label{eq:delayedCharacteristicEquation}
				\end{equation}
				where $H$ is given by Eq.~\eqref{eq:characteristicOL} and is the open-loop transfer function of the system without delays, i.e., for $\tau=0$. The roots of Eq.~\eqref{eq:delayedCharacteristicEquation} are the poles of the closed-loop system, and all of them must have a negative real part to guarantee the stability of this system \citep{Walton1987}. An inherent difficulty introduced by the presence of delays is that this characteristic equation is now transcendental. For nonzero $\tau$, the system possesses an infinity of poles which cannot be found in closed-form.
				
				\subsubsection{Root loci}
				
					Eq.~\eqref{eq:delayedCharacteristicEquation} can be solved numerically using homotopy: from the known solution at $\tau=0$ (where the characteristic equation is a polynomial), the root locus can progressively be computed. At each step, starting from a known solution for a given $\tau$, $\tau$ is incremented by $\Delta \tau$ and Eq.~\eqref{eq:delayedCharacteristicEquation} is solved with MATLAB's routine \texttt{fsolve} using as the initial guess the solution for the previous value of $\tau$. The procedure is then repeated until $\tau$ reaches a prescribed final value.
				
\begin{figure*}[!ht]
					\begin{subfigure}{.45\textwidth}				
						\centering
						\includegraphics[width=\textwidth]{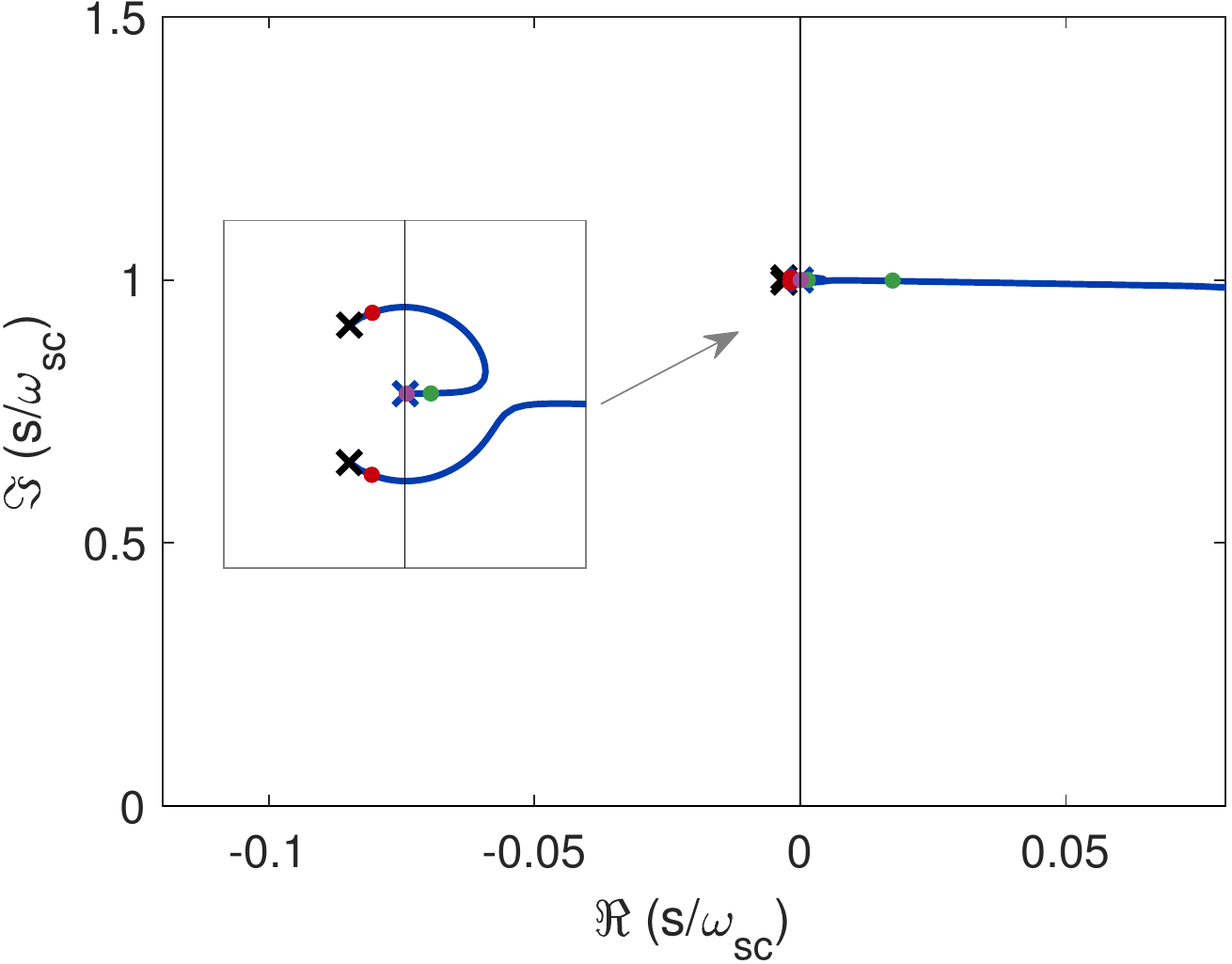}
						\caption{}
		  					\label{sfig:rootLocus_Kc1e-2}
					\end{subfigure}
					\hspace{0.05\textwidth}
					\begin{subfigure}{.45\textwidth}				
						\centering
						\includegraphics[width=\textwidth]{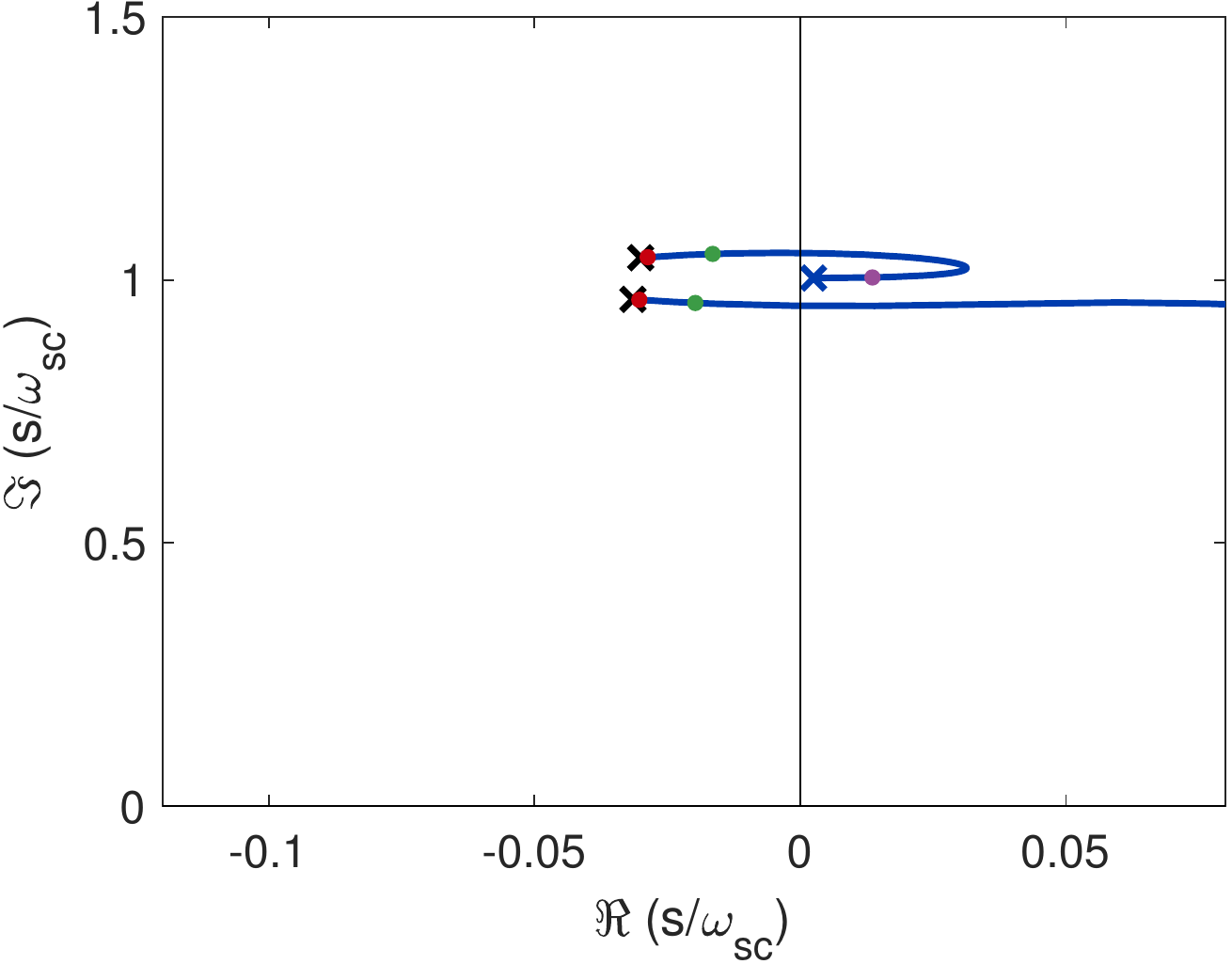}
						\caption{}
		  					\label{sfig:rootLocus_Kc1e-1}
					\end{subfigure}
					\caption{Root loci (parametrized by $\tau$) of the closed-loop system with delays ($\mathbf{\times}$: poles for $\tau = 0$, \textcolor{col2}{$\mathbf{\bullet}$}: $\tau = 0.01/\omega_{sc}$, \textcolor{col3}{$\mathbf{\bullet}$}: $\tau = 0.1/\omega_{sc}$, \textcolor{col4}{$\mathbf{\bullet}$}: $\tau = 1/\omega_{sc}$, \textcolor{col1}{$\mathbf{\times}$}: $\tau = \pi/\omega_{sc}$): $K_c = 0.01$~\subref{sfig:rootLocus_Kc1e-2} and  $K_c = 0.1$~\subref{sfig:rootLocus_Kc1e-1}.}
				\label{fig:rootLoci}
				\end{figure*}

			        Fig.~\ref{fig:rootLoci} shows root loci of the controlled system with delays for two values of $K_c$. The maximum value for $\tau$ is the maximum sampling period satisfying the Nyquist-Shannon sampling theorem if the system was forced at its resonant frequency $\omega_{sc}$, $\pi/\omega_{sc}$. %Only the poles with positive imaginary part are shown, but their complex conjugate counterpart also satisfy the characteristic equation.
				
				In all of these cases, the poles initially move to the right of the complex plane with increasing delays, and for large enough $\tau$ the highest-frequency poles cross the imaginary axis, which makes the closed-loop system unstable. As expected, the value of $\tau$ for which this instability occurs grows with $K_c$. %Intuitively, a controlled system with a higher $K_c$ will have poles which are further in the left half of the complex plane; a larger delay effect will thus be required to bring them to the right half.
				
				Fig.~\ref{fig:rootLoci} does not feature all the poles of the delayed system, except for $\tau=0$. As soon as $\tau>0$, a countable infinite set of poles emanate from $-\infty$, but these poles are not causing instabilities, unlike those featured in Fig.~\ref{fig:rootLoci}.
				
				\subsubsection{Critical delays}				
				
				Of particular interest is the value of $\tau$ for which the poles of the closed-loop system cross the imaginary axis, signalling the onset of instability. An inconvenient feature of Eq.~\eqref{eq:delayedCharacteristicEquation} is that this value can only be obtained by solving a transcendental equation. However, the following approximation can be made at frequencies comparatively low to the sampling frequency:
				\begin{multline}
					H_{ZOH}(s;\tau) = \dfrac{1-e^{-\tau s}}{\tau s} = e^{-\frac{\tau s}{2}} \dfrac{e^{\frac{\tau s}{2}} - e^{-\frac{\tau s}{2}}}{\tau s}  = e^{-\frac{\tau s}{2}} \dfrac{\displaystyle\sum_{k=0}^{+\infty} \left(\dfrac{\tau s}{2}\right)^k - \sum_{k=0}^{+\infty} \left(-\dfrac{\tau s}{2}\right)^k}{\tau s} \\ 
					= e^{-\frac{\tau s}{2}} \sum_{k=0}^{+\infty}\left(\dfrac{\tau s}{2}\right)^{2k}\approx e^{-\frac{\tau s}{2}},
				\end{multline}								
				i.e., the ZOH is nearly equivalent to a pure delay of $\tau/2$. With a pure delay model, the method of \cite{Walton1987} can be used to compute the characteristics roots of interest. Eq.~\eqref{eq:delayedCharacteristicEquation} is rewritten as				
				\begin{equation}
					1 + H(s)e^{-\frac{\tau s}{2}} = 0,
					\label{eq:delayedCharacteristicEquation2}
				\end{equation}
				The time delay resulting in purely imaginary characteristic roots is noted $\tau_c$. At this delay, a pair of complex conjugate poles or a single real pole cross the imaginary axis, possibly changing the stability of the system. Thus, $s=j\omega_c$ and $s=-j\omega_c$ satisfy the characteristic equation
				\begin{equation}
					\left\{\begin{array}{l} 
						1 + H(j\omega_c)e^{-\frac{j\omega_c \tau_c}{2}} = 0 \\
						1 + H(-j\omega_c)e^{\frac{j\omega_c \tau_c}{2}} = 0
					\end{array}\right. .
					\label{eq:characteristicComplexRoots}
				\end{equation}
				Multiplication of these two equations yield
				\begin{equation}
					H(j\omega_c) H(-j\omega_c) = 1.
					\label{eq:delayCondition1}
				\end{equation}
				This equation is a polynomial of $\omega_{sc}$; hence, there is a finite set of frequencies at which the poles of the closed-loop system cross the imaginary axis \citep{Walton1987}. The corresponding time delay $\tau_c$ can then be found using either line of Eq.~\eqref{eq:characteristicComplexRoots} as 
				\begin{equation}
					\tau_c = \dfrac{2}{\omega_c} \left[\angle -H(j\omega_c) + 2k\pi \right], \qquad k\in \mathbb{Z},
					\label{eq:delayCondition2}
				\end{equation}
				where $\angle$ is an operator giving the phase of a complex number.

			\subsubsection{Series approximations}
			
				Using Eq.~\eqref{eq:characteristicOL}, it can be shown that Eq.~\eqref{eq:delayCondition1} is a quartic polynomial of $\omega_c^2$. A closed-form solution can thus be obtained, but is impractically long. A more convenient form was obtained through Maclaurin series expansion in powers of $K_c$ of the analytical solution using Wolfram Mathematica. This provides an approximation of the critical frequencies. Among them, the one which corresponds to the minimum critical delay is given by
				\begin{equation}
						\omega_{c} = \omega_{sc} \left(1+K_c+\dfrac{5}{8}K_c^2 + \dfrac{73}{128}K_c^3 + O(K_c^4)\right).
				\end{equation}	
				Inserting this critical frequency into Eq.~\eqref{eq:delayCondition2} and expanding the result in power series of $K_c$ gives the corresponding critical delay
				\begin{equation}
						\tau_{c} = \dfrac{1}{\omega_{sc}}\left(\sqrt{6}\left(K_c-K_c^2\right) + \dfrac{19}{32}\sqrt{\dfrac{3}{2}}K_c^3 + O(K_c^4)\right).
					\label{eq:criticalDelaysSeries}
				\end{equation}
				$\tau_c$ corresponds to the largest admissible value of sampling time for a stable closed-loop system. It may also be noted that the first-order coefficient in $K_c$ obtained in Eq.~\eqref{eq:criticalDelaysSeries} for $\tau_{c}$ corresponds to that of the linearized value of $\tau$ leading to a non-passive delayed admittance (Eq.~\eqref{eq:tauNonPassive}).
				
				Eq.~\eqref{eq:criticalDelaysSeries} indicates that the critical sampling period is governed by the EEMCF of the system, whose value is typically small ($K_c\lesssim 0.1$). Hence, the associated critical sampling frequency may be orders of magnitude larger than the characteristic frequencies of the system. Although such a trend has been exhibited for other types of vibration absorbers before \citep{Olgac2000}, it is an important fact that needs to be accounted for when the DVA is used to emulate a passive shunt.
				
\begin{figure}[!ht]
					\centering
					\includegraphics[width=0.45\textwidth]{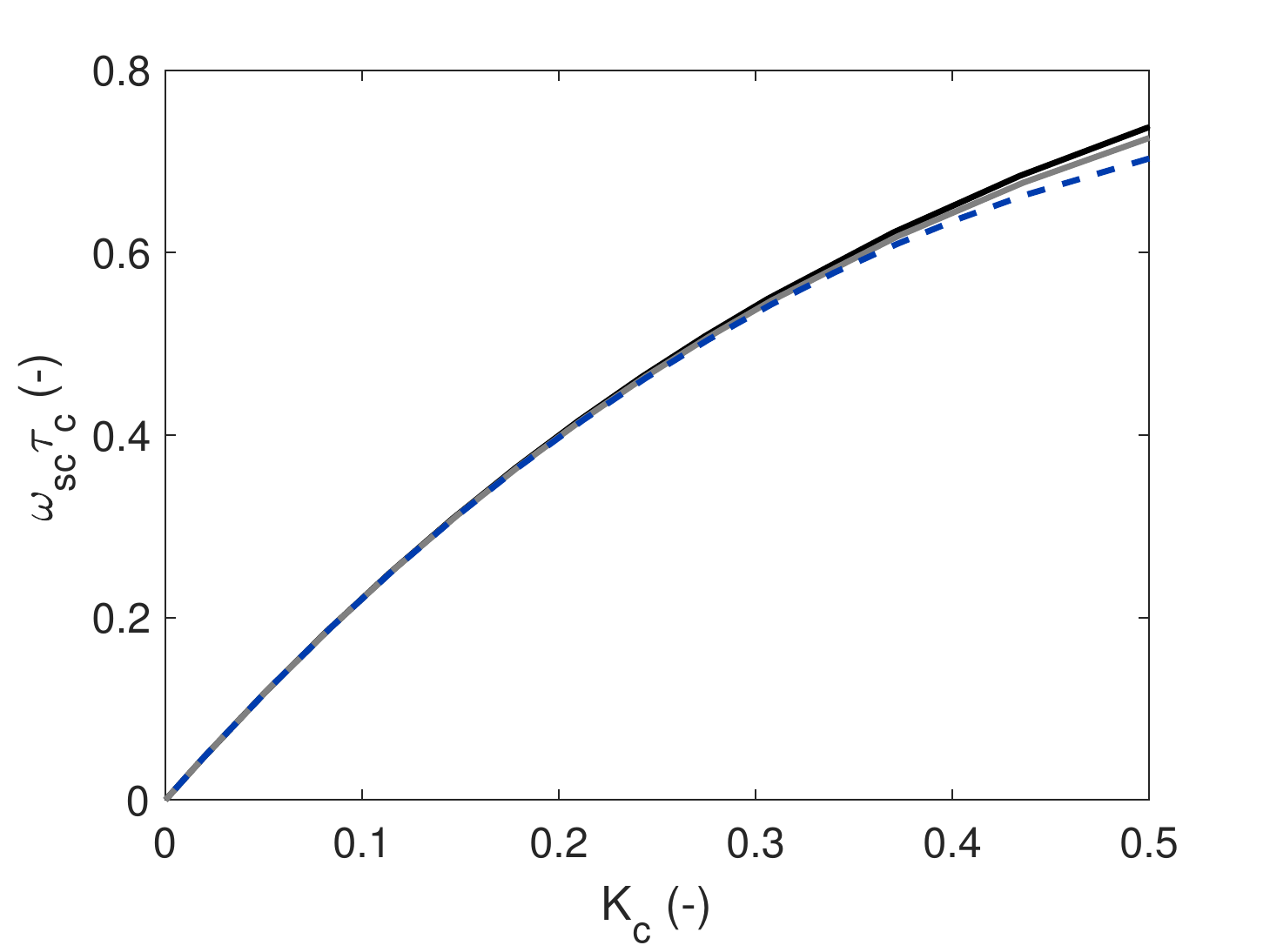}
					\caption{Critical delays $\tau_{c}$: ZOH model (\rule[0.2em]{1em}{0.2em}), pure delay model (\textcolor{gray1}{\rule[0.2em]{1em}{0.2em}}) and series approximation (\textcolor{col1}{\rule[0.2em]{0.4em}{0.2em}\nobreak\hspace{.2em}\rule[0.2em]{0.4em}{0.2em}}).}
					\label{fig:criticalDelay}
				\end{figure}

				The analytical approximations were compared with a direct numerical resolution of Eqs.~\eqref{eq:delayedCharacteristicEquation} and \eqref{eq:delayedCharacteristicEquation2} with $s=j\omega_c(K_c)$, solved with the \texttt{fsolve} routine from MATLAB using a homotopy on $K_c$. Fig.~\ref{fig:criticalDelay} compares the obtained critical delays. For small EEMCFs, the three models agree almost perfectly. Incidentally, this is also the range where the instabilities can be a compelling problem. %For large EEMCFs ($K_c \gtrsim 0.3$), the critical sampling frequency becomes a fraction of the Nyquist-Shannon frequency ($\tau = \pi/\omega_{sc}$), and thus delay-induced instabilities may not be the driving factor to set the sampling frequency.

		\subsubsection{FRF of the controlled system}
			
\begin{figure*}[!ht]
					\begin{subfigure}{.45\textwidth}				
						\centering
						\includegraphics[width=\textwidth]{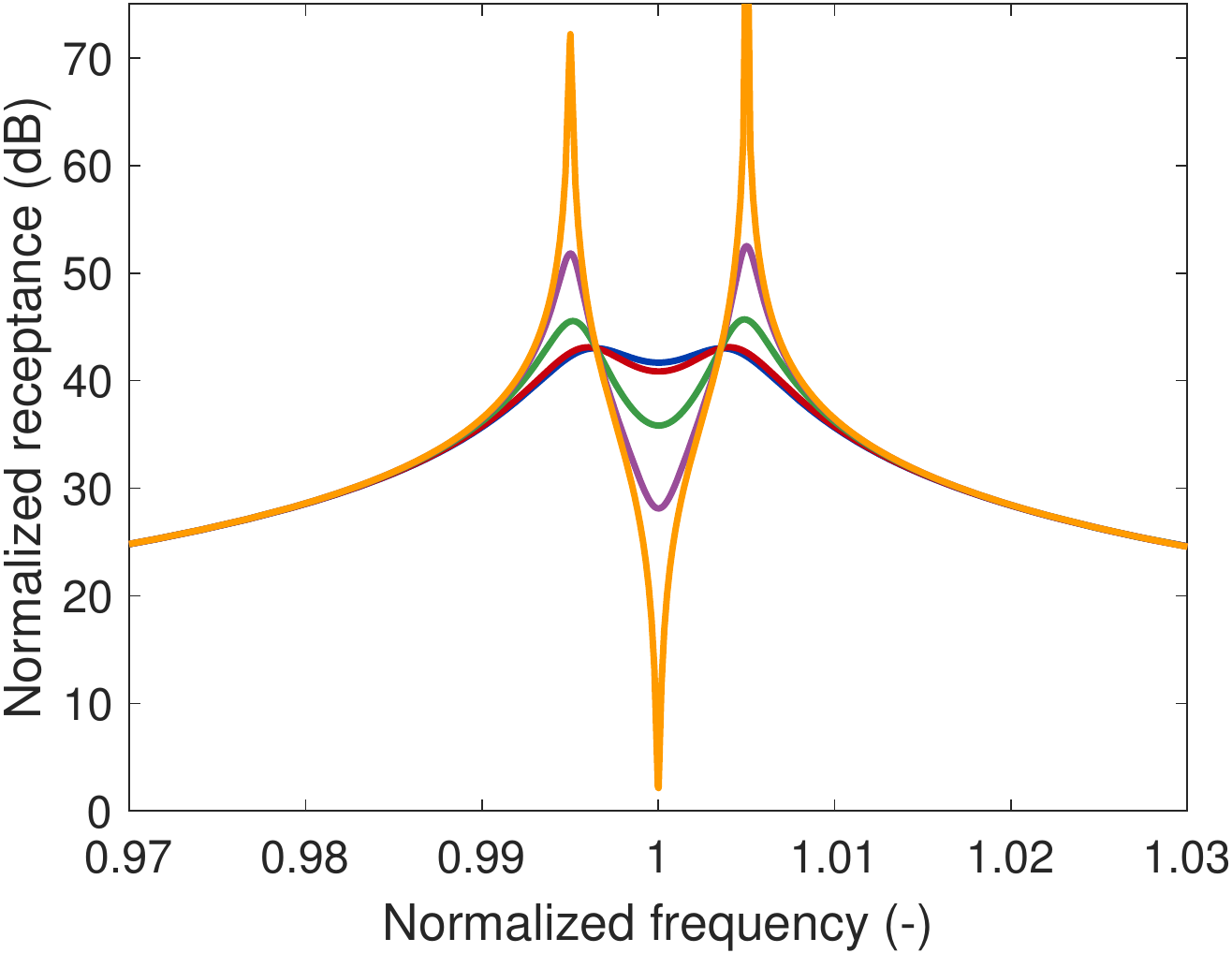}
						\caption{}
		  					\label{sfig:delayedFRF_Kc1e-2}
					\end{subfigure}
					\hspace{0.05\textwidth}
					\begin{subfigure}{.45\textwidth}				
						\centering
						\includegraphics[width=\textwidth]{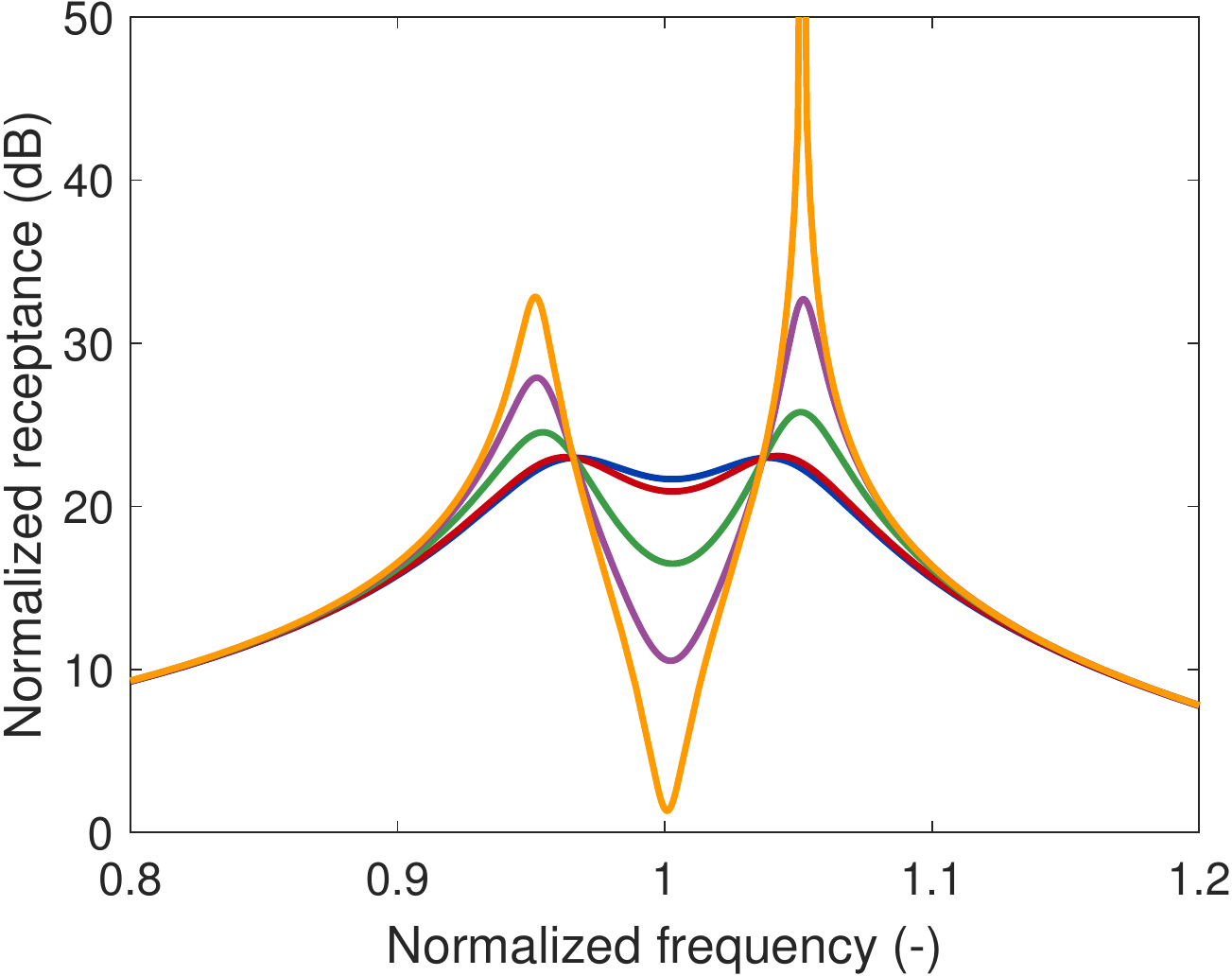}
						\caption{}
		  					\label{sfig:delayedFRF_Kc1e-1}
					\end{subfigure}
					\caption{FRF of the controlled system with a delayed admittance, $K_c=0.01$~\subref{sfig:delayedFRF_Kc1e-2} and $K_c=0.1$~\subref{sfig:delayedFRF_Kc1e-1}: $\tau = 0.01\tau_c$ (\textcolor{col1}{\rule[0.2em]{1em}{0.2em}}), $\tau = 0.1\tau_c$ (\textcolor{col2}{\rule[0.2em]{1em}{0.2em}}), $\tau = 0.5\tau_c$ (\textcolor{col3}{\rule[0.2em]{1em}{0.2em}}), $\tau = 0.8\tau_c$ (\textcolor{col4}{\rule[0.2em]{1em}{0.2em}}) and $\tau = \tau_c$ (\textcolor{col5}{\rule[0.2em]{1em}{0.2em}}).}
				\label{fig:delayedFRF}
				\end{figure*}

			Fig.~\ref{fig:delayedFRF} shows representative FRFs of the controlled system including the ZOH (using Eq.~\eqref{eq:zoh}). All the numerical FRFs and frequencies featured in this work are normalized with $k_{sc}$ and $\omega_{sc}$, respectively. Small sampling periods ($\tau \leq 0.1\tau_c$) have an imperceptible effect on the FRF compared to the nominal case. Conversely, a strong effect can be observed for large delays, especially on the rightmost peak whose amplitude grows with the delay. At $\tau=\tau_c$, the poles that lie on the imaginary axis create an undamped resonance in the FRF, signaling the onset of instability.

\begin{figure}[!ht]
					\centering
					\scalebox{0.7}{
					
        		\includegraphics[scale=1]{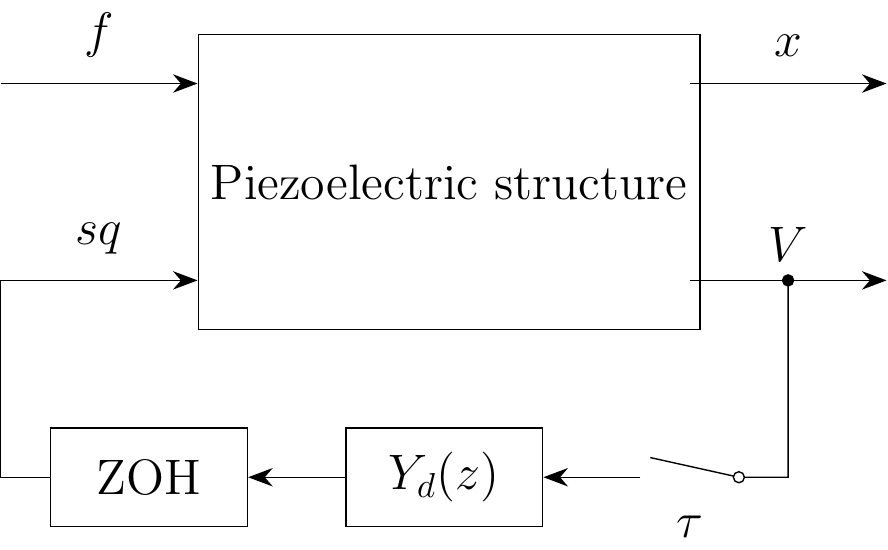}
					}
					\caption{Block diagram representation of the controlled system used for time simulations.}
					\label{fig:feedbackDelayedShuntSimulation}
				\end{figure}

\begin{figure*}[!ht]
					\begin{subfigure}{.45\textwidth}				
						\centering
						\includegraphics[width=\textwidth]{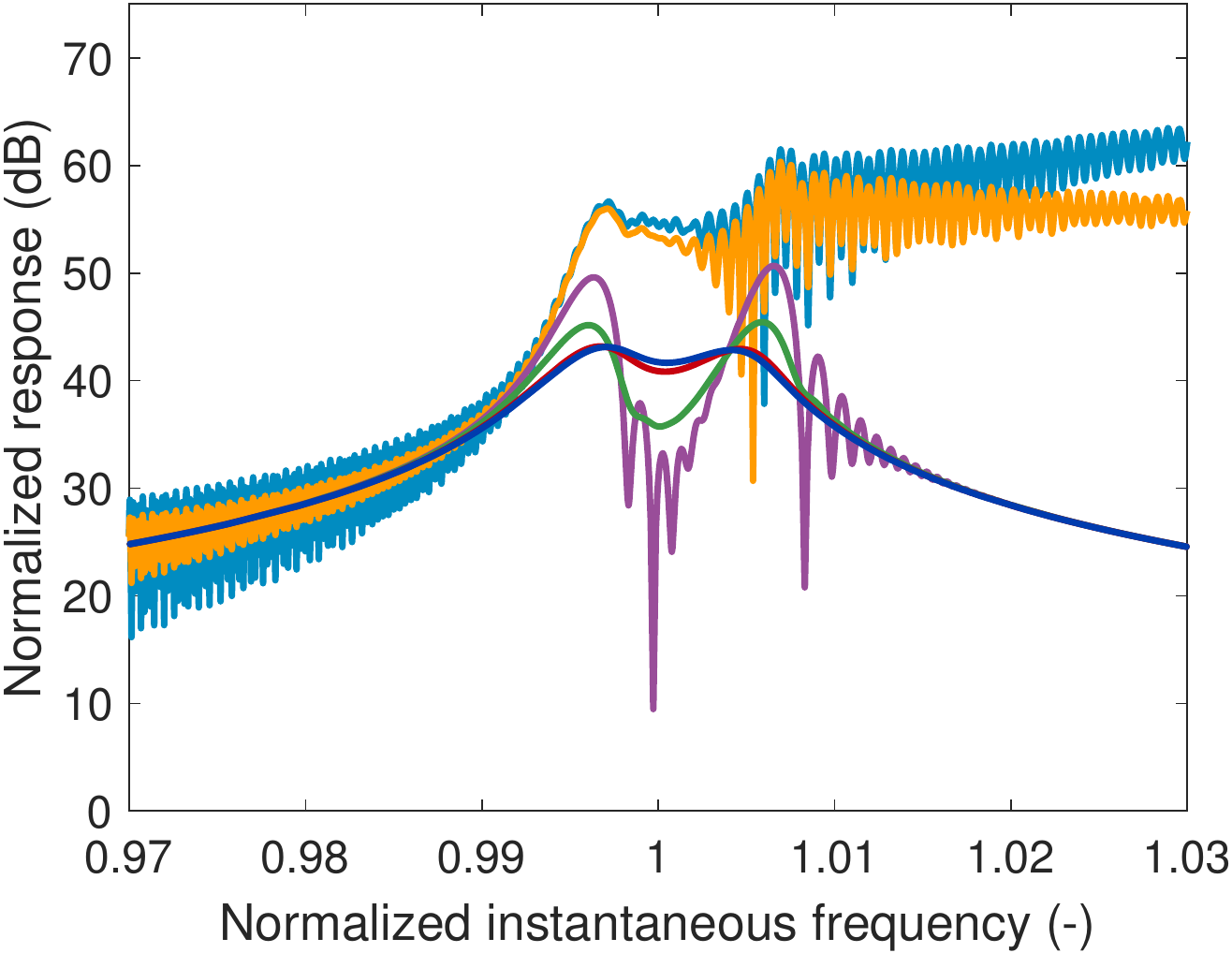}
						\caption{}
		  					\label{sfig:delayedTimeSimulation_Kc1e-2}
					\end{subfigure}
					\hspace{0.05\textwidth}
					\begin{subfigure}{.45\textwidth}				
						\centering
						\includegraphics[width=\textwidth]{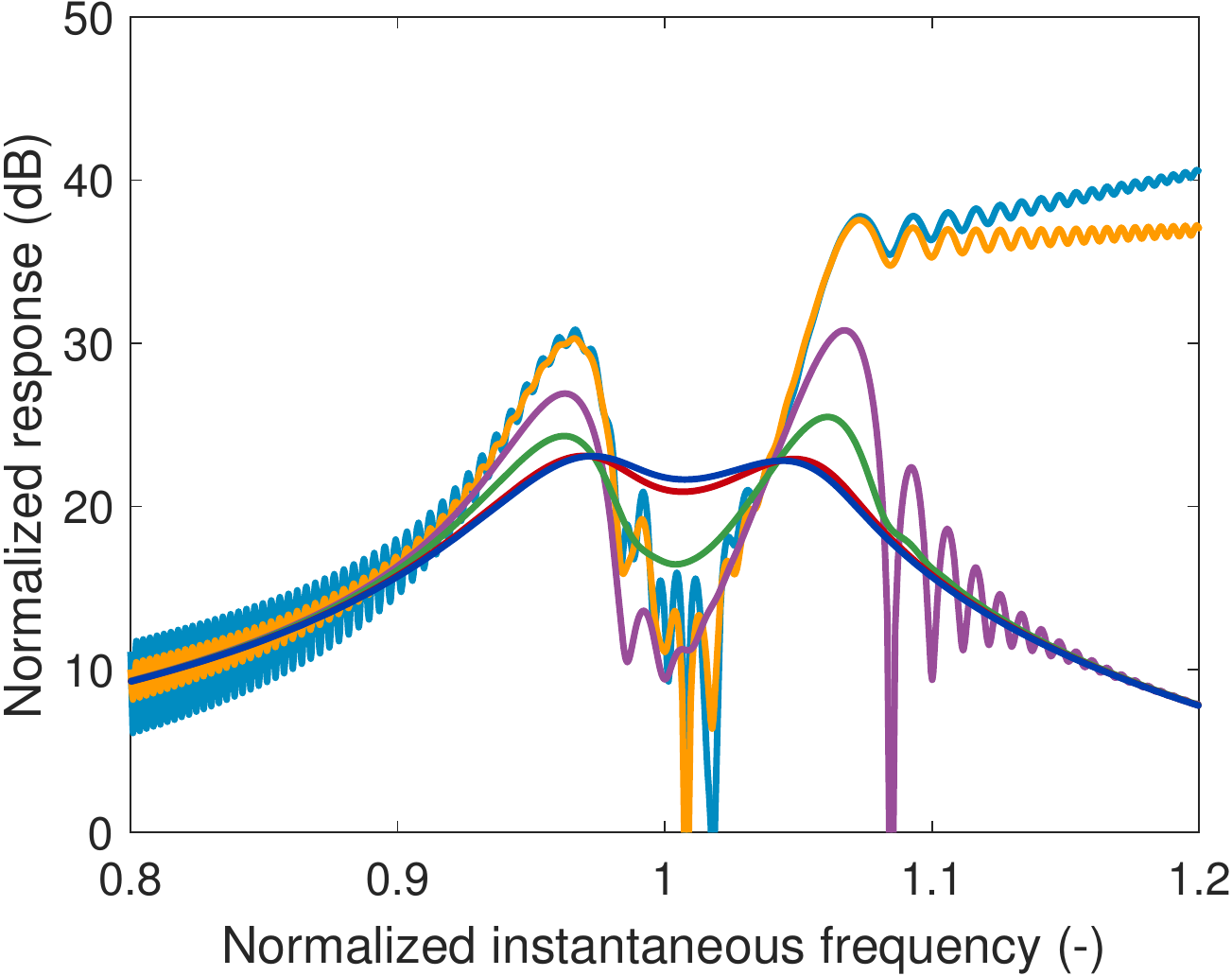}
						\caption{}
		  					\label{sfig:delayedTimeSimulation_Kc1e-1}
					\end{subfigure}
					\caption{Simulated envelope of the response of the controlled system with a delayed admittance to a unit-amplitude swept sine, $K_c=0.01$~\subref{sfig:delayedFRF_Kc1e-2} and $K_c=0.1$~\subref{sfig:delayedFRF_Kc1e-1}: $\tau = 0.01\tau_c$ (\textcolor{col1}{\rule[0.2em]{1em}{0.2em}}), $\tau = 0.1\tau_c$ (\textcolor{col2}{\rule[0.2em]{1em}{0.2em}}), $\tau = 0.5\tau_c$ (\textcolor{col3}{\rule[0.2em]{1em}{0.2em}}), $\tau = 0.8\tau_c$ (\textcolor{col4}{\rule[0.2em]{1em}{0.2em}}), $\tau = \tau_c$ (\textcolor{col5}{\rule[0.2em]{1em}{0.2em}}) and $\tau = 1.01\tau_c$ (\textcolor{col6}{\rule[0.2em]{1em}{0.2em}}).}
				\label{fig:simulatedDelayedFRF}
				\end{figure*}

				The analytical approximations were also verified by time simulations of the systems' responses to a unit-amplitude swept sine forcing under various sampling frequencies. The simulation of the system represented as a block diagram in Fig.~\ref{fig:feedbackDelayedShuntSimulation} was carried out with Simulink.	In addition to the sampling delay, this simulation accounts for the time-varying character of the system caused by sampling, as well as the effect of the discretization of the transfer function with Tustin's method. Fig.~\ref{fig:simulatedDelayedFRF} shows the envelopes of the systems' responses. The fact that the FRF is nearly not affected for $\tau \leq 0.1\tau_c$ is verified, and so is the progressive degradation, up to the onset of instability for $\tau\approx \tau_c$.
				
			According to the foregoing discussion, a rule of thumb is thus to choose the sampling period lower than or equal to one tenth of the critical delay. Besides, the sampling time must also be small enough so as to respect the Nyquist condition. Typical sampling frequencies of ten to thirty times the highest frequency of interest are often recommended \citep{Franklin1998}. The sampling time should therefore satisfy
			\begin{equation}
				\tau \leq \dfrac{1}{\omega_{sc}}\min\left\{\dfrac{2\pi}{30} , \dfrac{\sqrt{6}}{10}\left(K_c-K_c^2\right) + \dfrac{19}{320}\sqrt{\dfrac{3}{2}}K_c^3 \right\}.
				\label{eq:delays.maxSamplingPeriod}
			\end{equation}

	\section{Stabilization of delay-induced instabilities} 
	\label{sec:stabilization}
		\subsection{Discussion} 
		
			The delay-induced instabilities are clearly defeating the purpose of the DVA and should therefore be avoided. If the closed-loop system is prone to these instabilities, there are two possible options:
			\begin{enumerate}
				\item Choose a high enough sampling frequency.
				\item Modify the implemented admittance in anticipation of the delays.
			\end{enumerate}
			
			The first option is the most obvious and straightforward, but not always the most desirable one for two main reasons. 
			
			One reason is that a given digital unit's power consumption is a growing function of its clock frequency, which must be high enough to handle data at a given sampling frequency. The power consumption of the MCU can be estimated by \citep{Cardoso2017}
    		\begin{equation}
    			P_{\text{MCU}} = P_{\text{MCU,Static}}+P_{\text{MCU,Dynamic}} = V_{CC,MCU}I_{CC,MCU} + \beta_{MCU} C_L V_{CC,MCU}^2 f_{CPU}
    			\label{eq:digital.CPUPower}
    		\end{equation}
    		where $V_{CC,MCU}$ is the supply voltage, $I_{CC,MCU}$ the quiescent current, $\beta_{MCU}$ is the activity factor, $C_L$ is the load capacitance and $f_{CPU}$ is the clock frequency at which the digital unit is operating. Increasing the sampling frequency will increase $\beta_{MCU}$ and/or $f_{CPU}$, leading to a higher power consumption. Moreover, if $f_{CPU}$ is increased, $V_{CC,MCU}$ will also have to be increased, which leads to an actual power consumption proportional to $f_{CPU}^{3}$ \citep{Cardoso2017}.
			
			The second reason is that the required sampling frequency to make the delays effect negligible or let alone to have a stable closed-loop system may be very large. This would require high-frequency specialized equipments, whose cost may become prohibitively large.
		
		\subsection{Stabilization procedure} 		
		
			The principles of the proposed stabilization procedure are very similar to a pole placement approach: it is sought to place the poles of a modified delayed system as close as possible to those of the nominal system by modifying the parameters of the shunt admittance.

			\subsubsection{Pole placement via transfer function modification}
				
				The admittance of a shunt can be expressed as
				\begin{equation}
					Y_s(s) = \dfrac{\sum_{m=0}^M b_m s^m }{\sum_{n=0}^N a_n s^n}	
				\end{equation}							
				According to Eq.~\eqref{eq:characteristicOL2}, the poles of the nominal closed-loop system $p_k$ ($k=1,\cdots,K$) satisfy
				\begin{equation}
					1 - \dfrac{V(p_k)}{p_k q(p_k)} Y_s(p_k) = 0.
					\label{eq:stabilization.nominal}
				\end{equation}								
				In order to anticipate the delays, a modified admittance is introduced as
				\begin{equation}
					\widetilde{Y_s}(s) = \dfrac{\sum_{m=0}^M b_m (1+\delta_{b_m}) s^m }{\sum_{n=0}^N a_n (1+\delta_{a_n})s^n},
				\end{equation}
				where $\delta_{a_n}$ and $\delta_{b_m}$ are modification factors and are unknown for now. The poles of the modified delayed closed-loop system would be the solutions of Eq.~\eqref{eq:delayedCharacteristicEquation}:
				\begin{equation}
					1 - \dfrac{V(s)}{s q(s)}\dfrac{1-e^{-\tau s}}{\tau s} \widetilde{Y_s}(s) = 0.
					\label{eq:stabilization.modifiedDelayed}
				\end{equation}
				By comparing Eqs.~\eqref{eq:stabilization.nominal} and \eqref{eq:stabilization.modifiedDelayed}, in order for $p_k$ to be a pole of the modified delayed system, the modified delayed admittance must be equal to the nominal one at $s=p_k$:
				\begin{equation}
					\dfrac{1-e^{-\tau p_k}}{\tau p_k}\widetilde{Y_s}(p_k) = \dfrac{1-e^{-\tau p_k}}{\tau p_k}\dfrac{\sum_{m=0}^M b_m (1+\delta_{b_m}) p_k^m }{\sum_{n=0}^N a_n (1+\delta_{a_n})p_k^n}  =\dfrac{\sum_{m=0}^M b_m p_k^m }{\sum_{n=1}^N a_n p_k^n} = Y_s(p_k).
				\end{equation}							
				Rearranging this equation, the following relation is obtained
				\begin{equation}
					\dfrac{1-e^{-\tau p_k}}{\tau p_k}\dfrac{\sum_{m=0}^M b_m \delta_{b_m} p_k^m }{\sum_{m=0}^M b_m p_k^m } - \dfrac{\sum_{n=0}^N a_n \delta_{a_n} p_k^n }{\sum_{n=0}^N a_n p_k^n }  = 1 -\dfrac{1-e^{-\tau p_k}}{\tau p_k},
				\end{equation}							
				which, when imposed for $k=1,\cdots,K$, defines a linear system that can be put in a matrix form as
				\begin{equation}
					\begin{bmatrix}
						\mathbf{B} & \mathbf{A}
					\end{bmatrix}
					\begin{bmatrix}
						\delta_{b_0} \\ \vdots \\ \delta_{b_M} \\ \delta_{a_0} \\ \vdots \\ \delta_{a_N}
					\end{bmatrix} =
					\begin{bmatrix}
						1 -\dfrac{1-e^{-\tau p_1}}{\tau p_1} \\
						\vdots \\
						1 -\dfrac{1-e^{-\tau p_K}}{\tau p_K}
					\end{bmatrix},
					\label{eq:delaysSystem}
				\end{equation}	
				where
				\begin{equation}
				    \mathbf{B} = \begin{bmatrix}
						\dfrac{1-e^{-\tau p_1}}{\tau p_1} \dfrac{b_0}{\displaystyle \sum_{m=0}^M b_m p_1^m} & \cdots & \dfrac{1-e^{-\tau p_1}}{\tau p_1} \dfrac{b_M p_1^M}{\displaystyle \sum_{m=0}^M b_m p_1^m}\\
						\vdots & & \vdots  \\
						\dfrac{1-e^{-\tau p_K}}{\tau p_K} \dfrac{b_0}{\displaystyle \sum_{m=0}^M b_m p_K^m} & \cdots & \dfrac{1-e^{-\tau p_K}}{\tau p_K} \dfrac{b_M p_K^M}{\displaystyle \sum_{m=0}^M b_m p_K^m}
					\end{bmatrix}
				\end{equation}
				and
				\begin{equation}
				    \mathbf{A} = \begin{bmatrix}
						- \dfrac{a_0 }{\displaystyle \sum_{n=0}^N a_n p_1^n } & \cdots & - \dfrac{a_N p_1^N }{\displaystyle \sum_{n=0}^N a_n p_1^n } \\
						\vdots & & \vdots \\
						- \dfrac{a_0 }{\displaystyle \sum_{n=0}^N a_n p_K^n } & \cdots & - \dfrac{a_N p_K^N }{\displaystyle \sum_{n=0}^N a_n p_K^n }
					\end{bmatrix}.
				\end{equation}
				In short, Eq.~\eqref{eq:delaysSystem} can be rewritten
				\begin{equation}
					\mathbf{P}\bm{\delta} = \mathbf{d}.
					\label{eq:delays_coeffModSystem}
				\end{equation}
				This system has a trivial solution $\bm{\delta} = [-1,\cdots,-1]^T$. This makes all the coefficients of the modified admittance equal to zero, which clearly is not an acceptable solution. To resolve this, one of the modification factor can be imposed to an arbitrary value, for instance 0. For this particular choice, the column associated with this modification coefficient may simply be removed from $\mathbf{P}$. Thus, the number of unknowns is reduced to $M+N+1$. Since this number may not be equal to $K$, the system may not be square. To solve it, the pseudoinverse (denoted by a superscript $\dagger$) is used.
				\begin{equation}
					\bm{\delta} = \mathbf{P}^\dagger \mathbf{d}.
					\label{eq:delays_pseudoinverse}
				\end{equation}
				%The formulation presented herein differs slightly from that presented in~\cite{Raze2019}, in that the unknowns are modification factors rather than the coefficients of the modified transfer function. The latter generally have widely different scales, which can make the problem badly conditioned, and in case of an underdetermined system, some coefficient may undergo a very large relative modification. By contrast, the formulation used here allows to obtain a somewhat balanced relative modification of the coefficients. 
				It should be noted that the procedure only requires the knowledge of the sampling period $\tau$ in addition to what is already known for tuning the shunt. This parameter is set by the user and is thus well-known and well-controlled.
				
				\subsubsection{Numerical verification}
				
\begin{figure*}[!ht]
					\begin{subfigure}{.45\textwidth}				
						\centering
						\includegraphics[width=\textwidth]{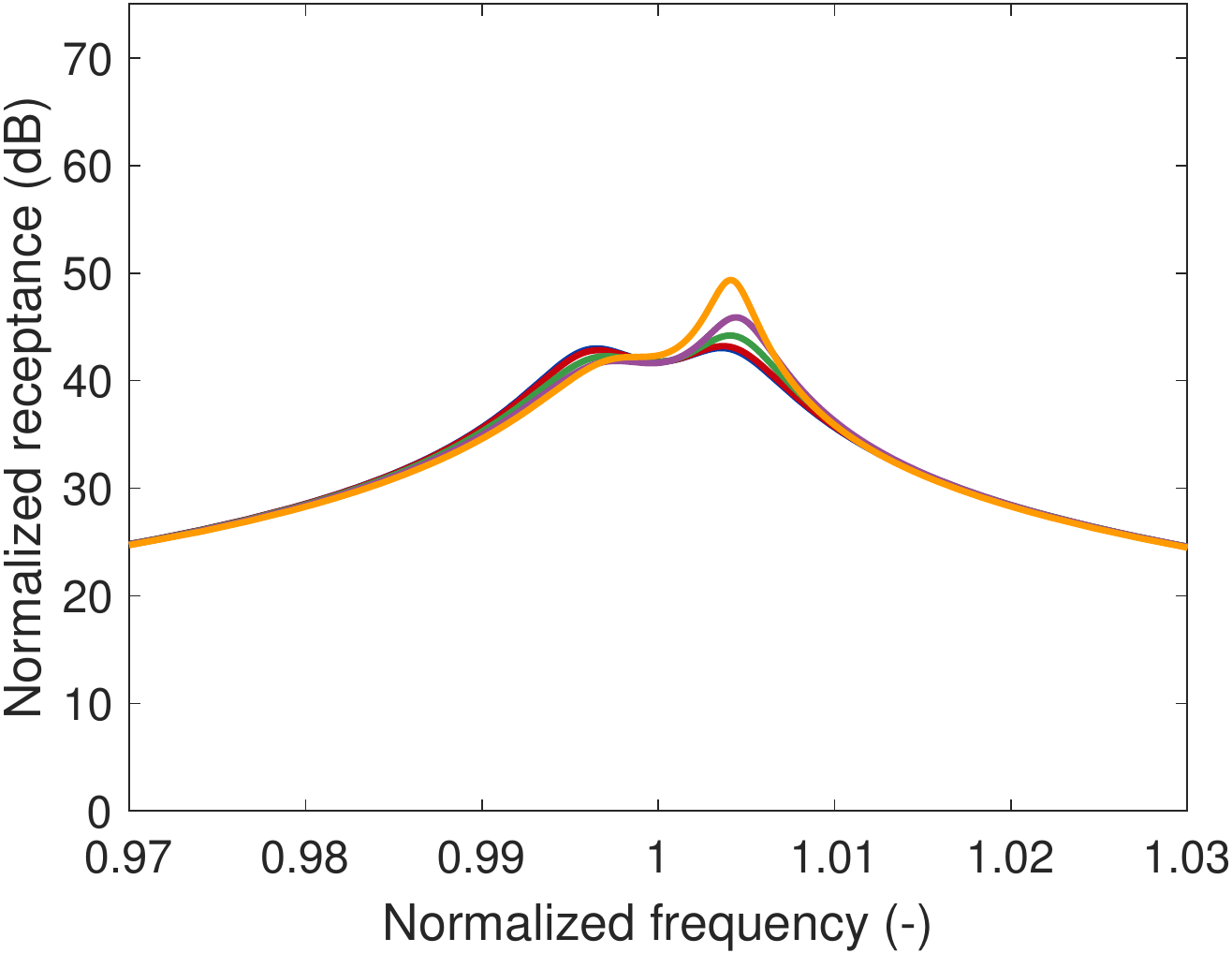}
						\caption{}
		  					\label{sfig:delayedFRF_Kc1e-2mod}
					\end{subfigure}
					\hspace{0.05\textwidth}
					\begin{subfigure}{.45\textwidth}				
						\centering
						\includegraphics[width=\textwidth]{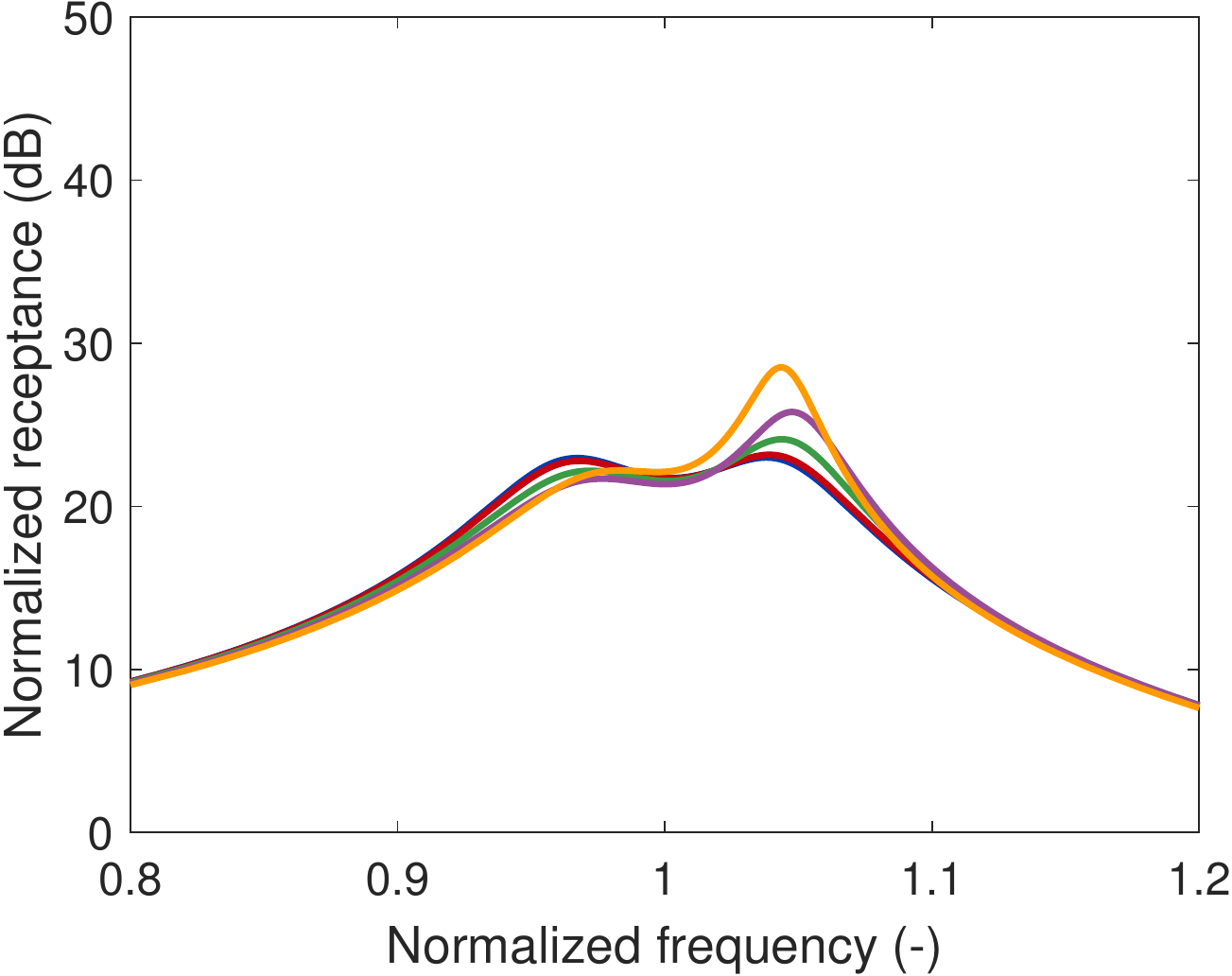}
						\caption{}
		  					\label{sfig:delayedFRF_Kc1e-1mod}
					\end{subfigure}
					\caption{FRF of the controlled system with a modified delayed admittance, $K_c=0.01$~\subref{sfig:delayedFRF_Kc1e-2mod} and $K_c=0.1$~\subref{sfig:delayedFRF_Kc1e-1mod}: $\tau = 0.01/\omega_{sc}$ (\textcolor{col1}{\rule[0.2em]{1em}{0.2em}}), $\tau = 0.1/\omega_{sc}$ (\textcolor{col2}{\rule[0.2em]{1em}{0.2em}}), $\tau = 0.5/\omega_{sc}$ (\textcolor{col3}{\rule[0.2em]{1em}{0.2em}}), $\tau = 1/\omega_{sc}$ (\textcolor{col4}{\rule[0.2em]{1em}{0.2em}}) and $\tau = \pi/\omega_{sc}$ (\textcolor{col5}{\rule[0.2em]{1em}{0.2em}}).}
				\label{fig:delayedFRFmod}
				\end{figure*}

Fig.~\ref{fig:delayedFRFmod} shows the results of the stabilization procedure on the FRF of the controlled system. The maximum sampling period $\pi/\omega_{sc}$, is about 130 and 14 times greater than $\tau_c$ for $K_c = 0.01$ and $K_c = 0.1$, respectively. A remarkable feature is that two FRFs for a different EEMCF but with an identical sampling frequency look similar, unlike the unmodified case. Therefore, with this modification, the EEMCF no longer appears to play a role in the delay-induced degradation of the vibration reduction. The FRFs for $\tau > 0.1/\omega_{sc}$ do not exhibit as good performance as the others, but it is not advised to choose such a low sampling frequency anyway \citep{Franklin1998}.
				
\begin{figure*}[!ht]
					\begin{subfigure}{.45\textwidth}				
						\centering
						\includegraphics[width=\textwidth]{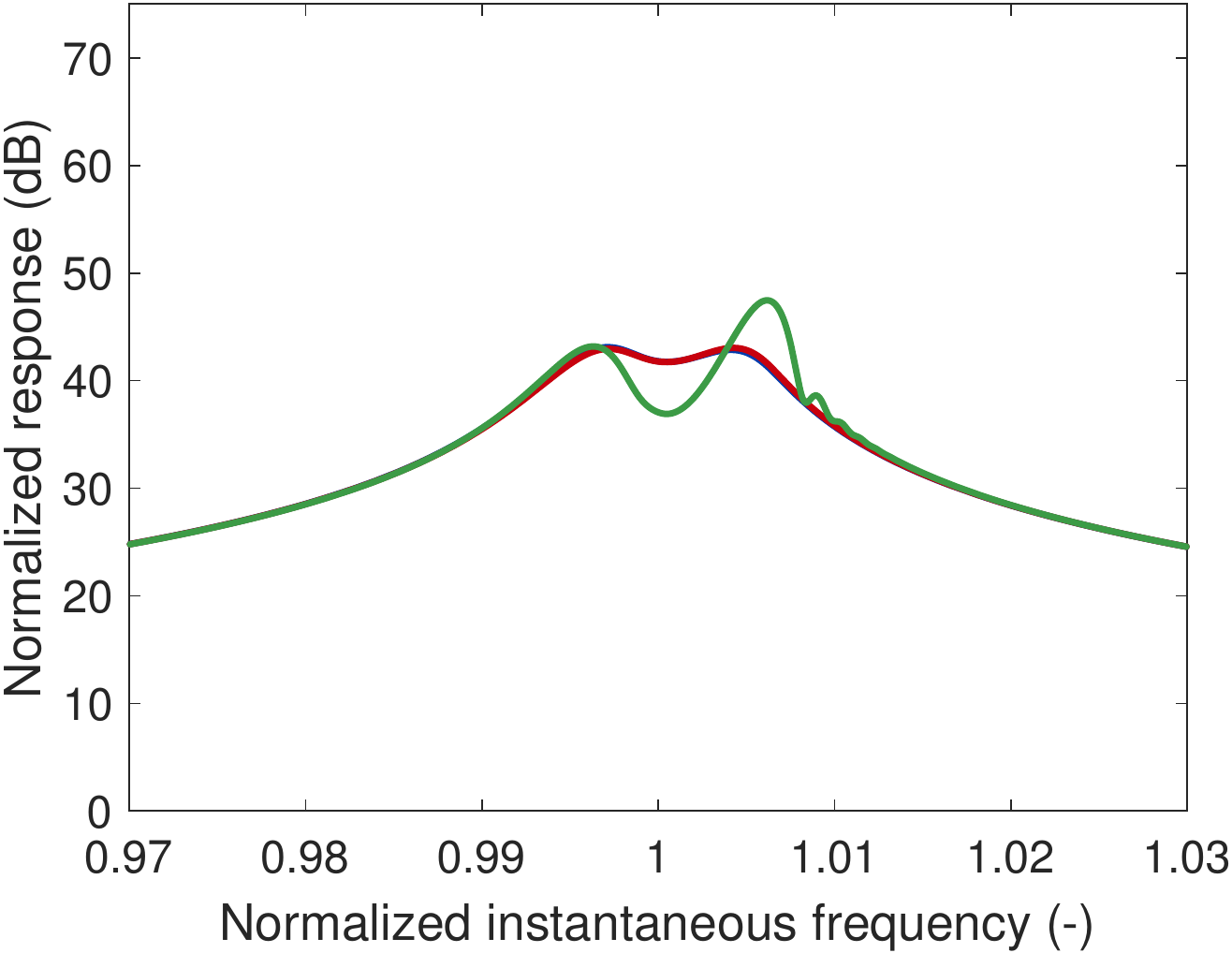}
						\caption{}
		  					\label{sfig:delayedFRF_Kc1e-2modSim}
					\end{subfigure}
					\hspace{0.05\textwidth}
					\begin{subfigure}{.45\textwidth}				
						\centering
						\includegraphics[width=\textwidth]{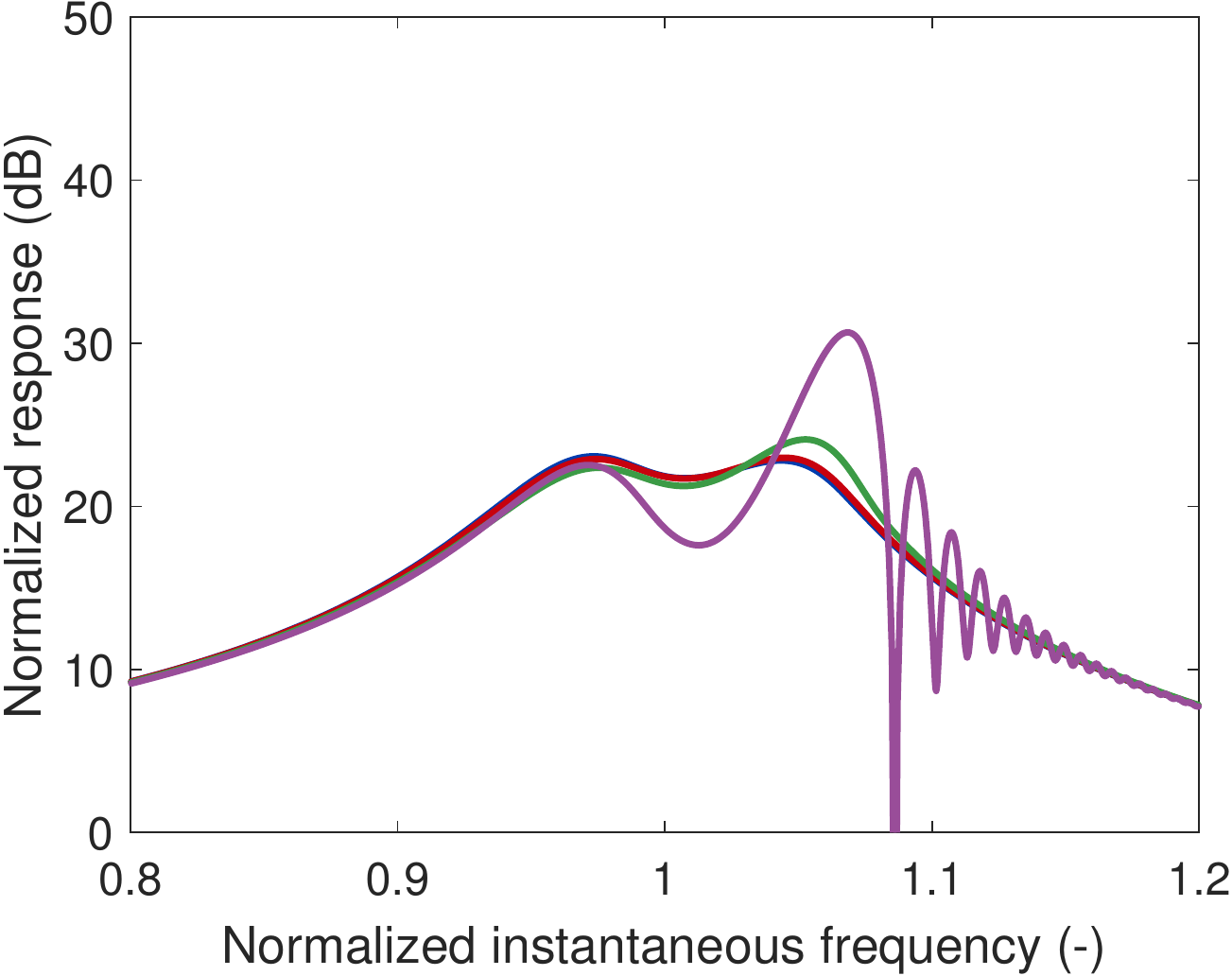}
						\caption{}
		  					\label{sfig:delayedFRF_Kc1e-1modSim}
					\end{subfigure}
					\caption{Simulated envelope of the response of the controlled system with a delayed, modified admittance to a unit-amplitude swept sine, $K_c=0.01$~\subref{sfig:delayedFRF_Kc1e-2modSim} and $K_c=0.1$~\subref{sfig:delayedFRF_Kc1e-1modSim}: $\tau = 0.01/\omega_{sc}$ (\textcolor{col1}{\rule[0.2em]{1em}{0.2em}}), $\tau = 0.1/\omega_{sc}$ (\textcolor{col2}{\rule[0.2em]{1em}{0.2em}}), $\tau = 0.5/\omega_{sc}$ (\textcolor{col3}{\rule[0.2em]{1em}{0.2em}}) and $\tau = 1/\omega_{sc}$ (\textcolor{col4}{\rule[0.2em]{1em}{0.2em}}).}
				\label{fig:delayedFRFmodSimulation}
				\end{figure*}

				A second verification was made with the time simulation of the system featured in Fig.~\ref{fig:feedbackDelayedShuntSimulation} using the modified admittance parameters. By comparing Figs.~\ref{fig:delayedFRFmod} and \ref{fig:delayedFRFmodSimulation}, it can be observed that both models agree well for $\tau \leq 0.1/\omega_{sc}$. However, discrepancies appear above this limit. In particular, the systems are unstable for sampling periods equal to or larger than $\tau = 1/\omega_{sc}$ and $\tau = \pi/\omega_{sc}$ for $K_c=0.01$ and $K_c = 0.1$, respectively. This can be attributed to the time-variant characteristics of sampling which were neglected in the analysis, as well as the frequency warping due to Tustin's transform. In practice, it is not advised to choose a sampling frequency smaller than thirty times the highest frequency of interest \citep{Franklin1998}. Eq.~\eqref{eq:delays.maxSamplingPeriod} therefore becomes
				\begin{equation}
					\tau \leq \dfrac{2\pi}{30\omega_{sc}}
				\end{equation}
				to ensure the stability of the closed-loop system with a modified admittance with some margin.
				
		\section{Experimental validation}		
		\label{sec:experiments}
		
			\begin{figure*}[!ht]
    			\begin{subfigure}{\textwidth}
    		    \centering
    				\scalebox{0.7}{
        		\includegraphics[scale=1]{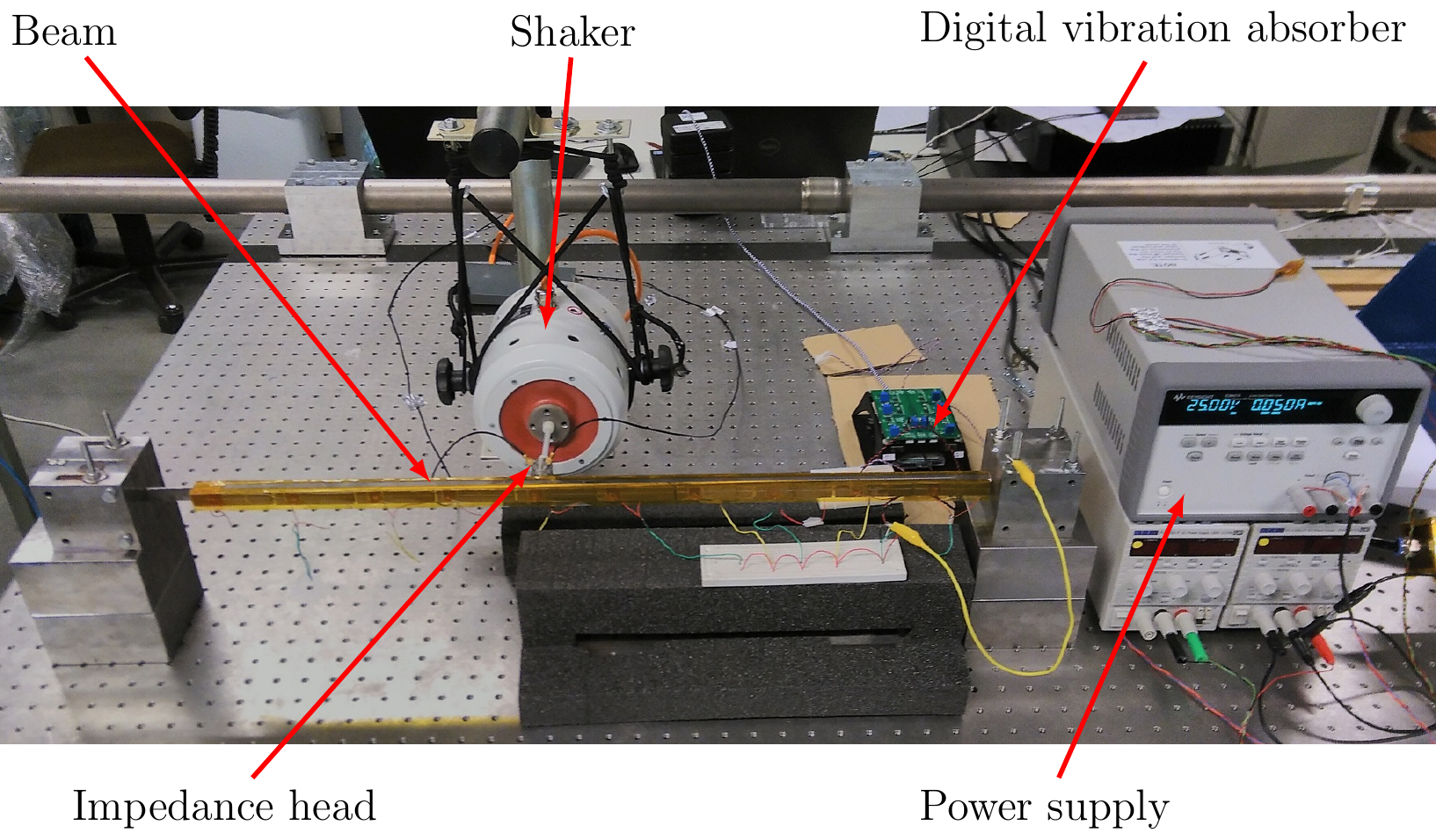}
    				}
    				\caption{}
    				\label{sfig:setupPicture}
    			\end{subfigure}
    			
    		    \vspace{1em}
    			\begin{subfigure}{\textwidth}
    				\centering
    				\scalebox{0.7}{
        		\includegraphics[scale=1]{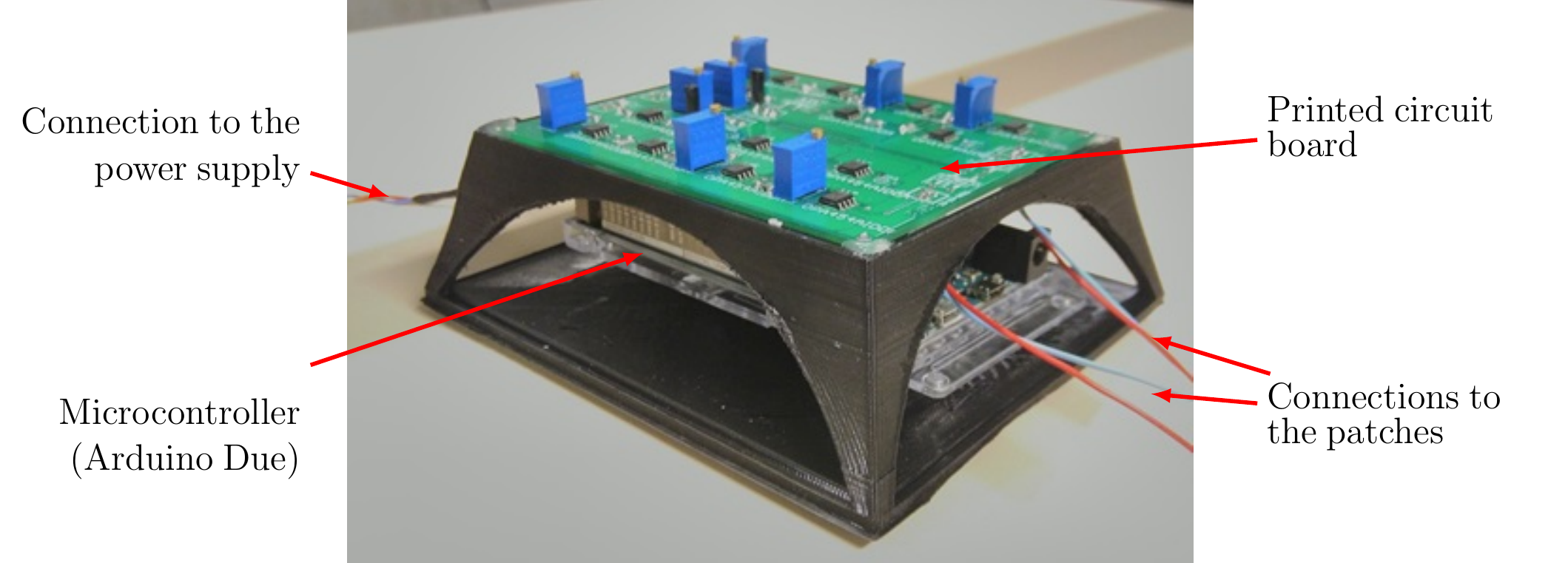}
    				}
    				\caption{}
    				\label{sfig:digital.DAPicture1}
    			\end{subfigure}
    		    
    		    \caption{Picture of the experimental setup~\subref{sfig:setupPicture} and close-up on the DVA~\subref{sfig:digital.DAPicture1}.}
    		    \label{fig:setup}
    		\end{figure*}				
			
			The action of the DVA was experimentally demonstrated on a clamped-free beam with a clamped thin lamina attached to its free end, as shown in Fig.~\ref{fig:setup}\subref{sfig:setupPicture}.  It was excited at middle span by an electrodynamic shaker (TIRA TV 51075). An impedance head (DYTRAN 5860B) was used to measure the force applied to the structure and the corresponding acceleration. The measurements were recorded by an acquisition system (LMS Scadas Mobile). The beam is covered over its whole length with ten cells, composed of pairs of stacks of two PSI-5A4E piezoelectric patches, each pair being placed on either side of the beam. The five cells closest to the clamped end were connected in parallel and used for mitigating the resonant vibrations around the first beam mode, whereas the five other cells were left in open circuit. More details about the experimental setup can be found in~\cite{Lossouarn2018}.
			
    		 To identify the system, the FRFs of the beam whith short-circuited and open-circuited patches were measured. Fitting these FRFs gave an estimation of the short- and open-circuit resonance frequencies. The piezoelectric capacitance was then measured with a multimeter (FLUKE 177). From these parameters, the optimal inductance and resistance of a series RL shunt were computed using Eqs.~\eqref{eq:SeriesLopt} and \eqref{eq:SeriesRopt}, respectively. All these parameters are reported in Table~\ref{tab:experimentalSetupData}.
    		 
			\begin{table*}[!ht]
				\small\sf\centering
				\caption{Parameters of the experimental setup.}
				\begin{tabular}{ccccccc}
				\hline
					Parameter & $f_{sc}$ & $f_{oc}$ & $K_c$ & $C_p^\varepsilon$ &  $R$ & $L$\\
					\hline
					Value & 31.08Hz & 31.29Hz & 0.116 & 245nF &  2,961$\Omega$ & 105.7H \\
					\hline
				\end{tabular}
				\label{tab:experimentalSetupData}
			\end{table*}
			
			The DVA shown in  Fig.~\ref{fig:setup}\subref{sfig:digital.DAPicture1} was powered with $\pm 25V$, and the MCU was programmed in order to mimic the admittance of the series RL shunt. The resistance of the current injector $R_i = 268.3\Omega$ was measured with a multimeter.

\begin{figure*}[!ht]
					\begin{subfigure}{.45\textwidth}				
						\centering
						\includegraphics[width=\textwidth]{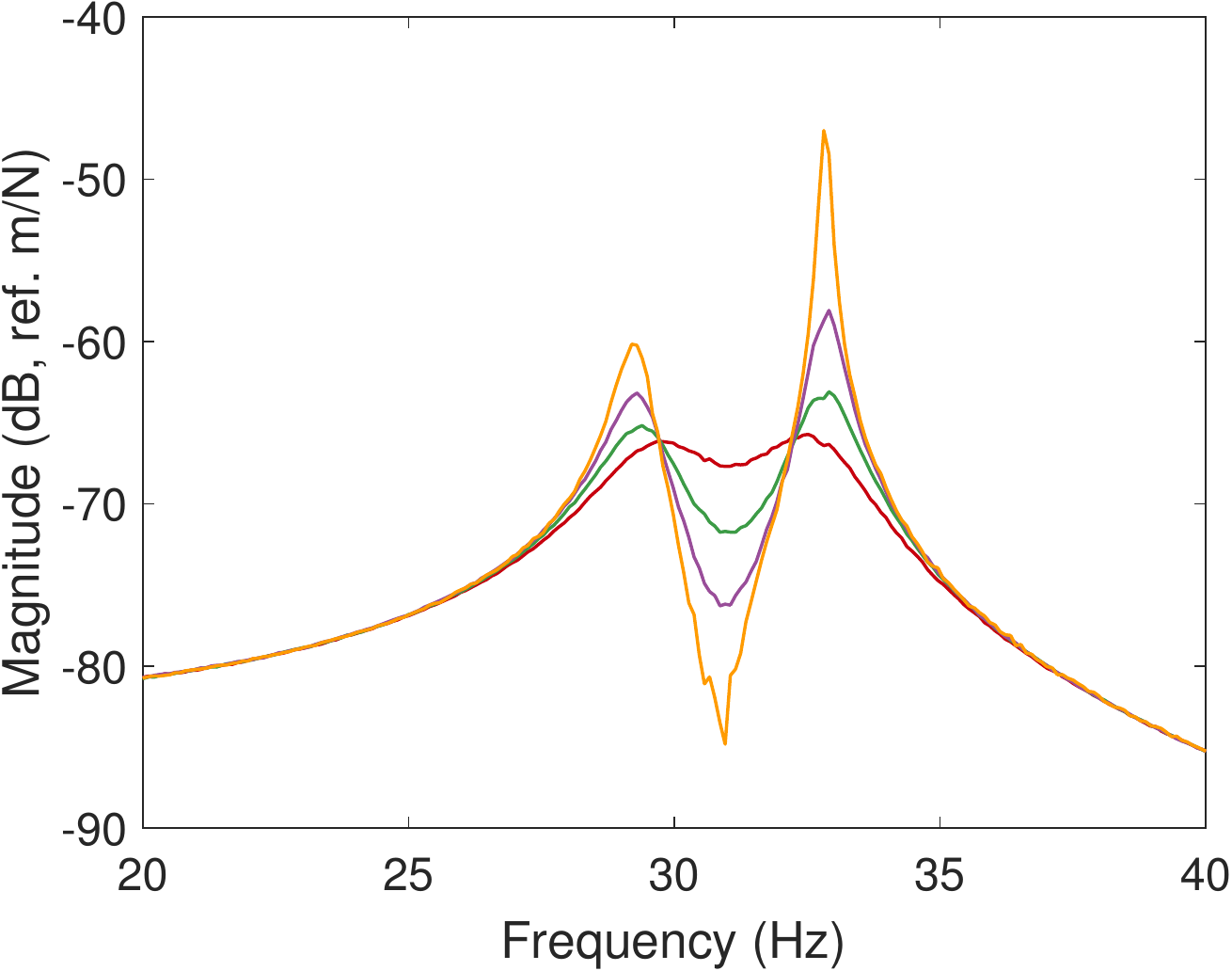}
						\caption{}
		  					\label{sfig:delayedFRFExp}
					\end{subfigure}
					\hspace{0.05\textwidth}
					\begin{subfigure}{.45\textwidth}				
						\centering
						\includegraphics[width=\textwidth]{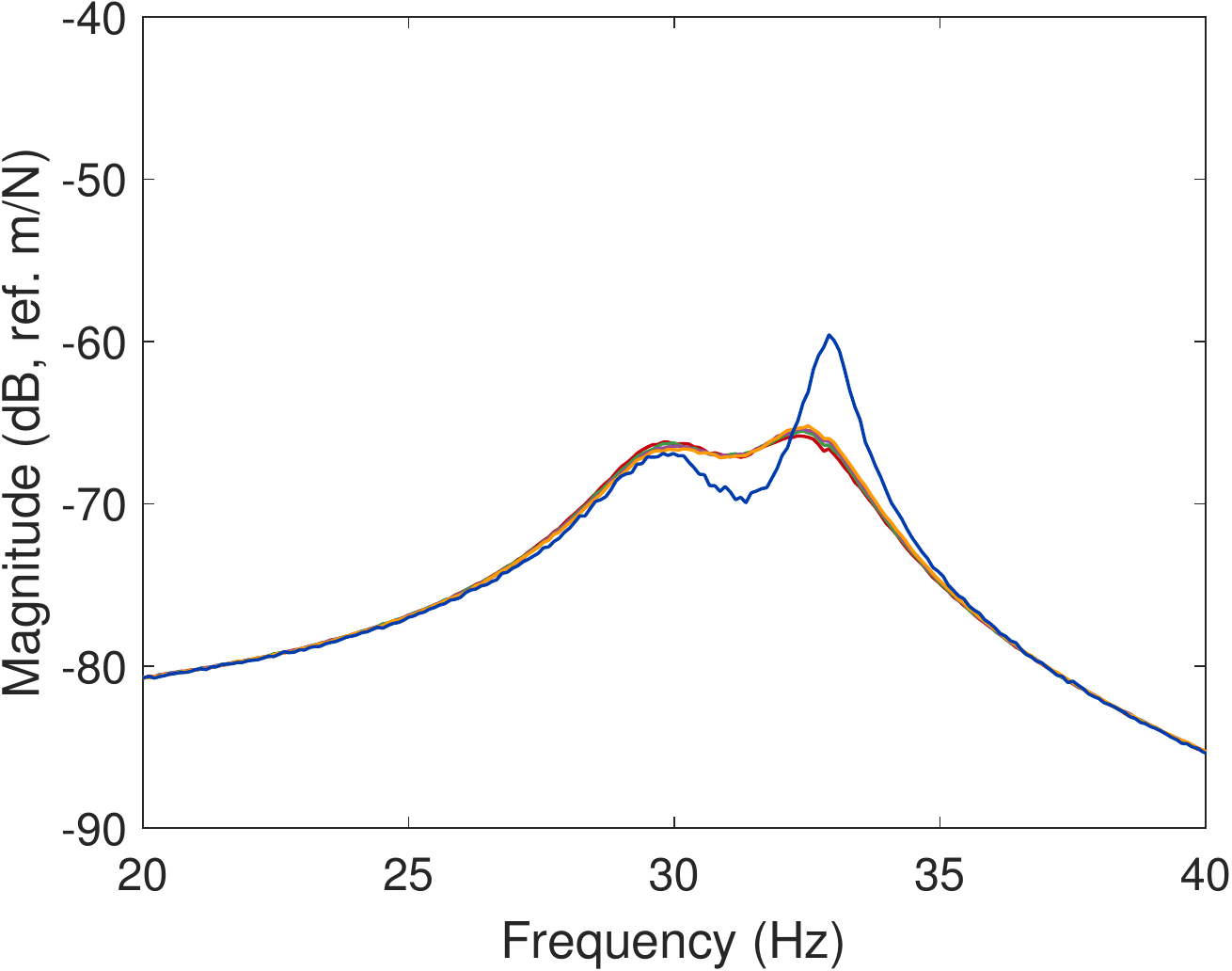}
						\caption{}
		  					\label{sfig:delayedFRFExpMod}
					\end{subfigure}
					\caption{Experimental FRF of the beam ($K_c=0.116$) with an unmodified~\subref{sfig:delayedFRFExp} and a modified~\subref{sfig:delayedFRFExpMod} admittance: $\tau = 10^{-4}$s$\approx 0.1\tau_c$ (\textcolor{col2}{\rule[0.2em]{1em}{0.2em}}), $\tau = 6.5\times 10^{-4}$s$=0.5\tau_c$ (\textcolor{col3}{\rule[0.2em]{1em}{0.2em}}), $\tau = 10^{-3}$s$\approx 0.8\tau_c$ (\textcolor{col4}{\rule[0.2em]{1em}{0.2em}}), $\tau = 1.3\times 10^{-3}$s$=\tau_c$ (\textcolor{col5}{\rule[0.2em]{1em}{0.2em}}) and $\tau = 5\times 10^{-3}$s$ \approx 1/\omega_{sc}$ (\textcolor{col1}{\rule[0.2em]{1em}{0.2em}}) .}
				\label{fig:delayedFRFExp}
				\end{figure*}

			To experimentally validate the developments about delay-induced instabilities, FRFs were measured under progressively decreasing sampling frequencies. As testified by Fig.~\ref{fig:delayedFRFExp}\subref{sfig:delayedFRFExp}, the destabilization effect of the sampling frequency is clearly observable. The results featured in this figure are close to those of Fig.~\ref{fig:delayedFRF}\subref{sfig:delayedFRF_Kc1e-1} (the coupling factor of the experimental setup is 0.116, which is close to the EEMCF of 0.1 used therein), which validates the model used to describe sampling delays. From Eq.~\eqref{eq:criticalDelaysSeries}, the stability limit of the unmodified system should theoretically be reached at $\tau = 1.3\times 10^{-3}$s. The experimental system is still stable but very lightly damped. This small discrepancy can be explained by the presence of structural damping in the host, as well as by experimental uncertainties.
			
			 The stabilization procedure recovers the performance of a case without delays, as shown in Fig.~\ref{fig:delayedFRFExp}\subref{sfig:delayedFRFExpMod}. Namely, all the curves are virtually superimposed up to $\tau = \tau_c$, which validates the proposed stabilization method. Fourfold a sampling period leads to a system with modified admittance where the effects of sampling are observable, more than in the numerical model featured in Fig.~\ref{fig:delayedFRFmod}\subref{sfig:delayedFRF_Kc1e-1mod}, but similarly to the time simulation in Fig.~\ref{fig:delayedFRFmodSimulation}\subref{sfig:delayedFRF_Kc1e-1modSim}. Nevertheless, a case with such a high sampling period when the admittance is unmodified is not disclosed here, as it leads to an unstable closed-loop system. 
			
    \section{Conclusion}
    \label{sec:conclusion}
    
        A DVA used for piezoelectric shunt damping is an attractive solution but it may be hindered by delay-induced instabilities incurred by the sampling procedure in the digital unit. After reviewing the basics of piezoelectric shunt damping with a DVA, this work used concepts from feedback control theory to highlight the small phase margin exhibited by piezoelectric systems with small EEMCFs. This makes them susceptible to delay-induced instabilities when a DVA is used, despite the passive character of the control law. An approximate explicit relation was derived between the maximum sampling period guaranteeing stability and the EEMCF, and it was shown that this period tends to zero concurrently with the EEMCF.
        
        Since piezoelectric structures typically exhibit small EEMCFs, the maximum sampling period may be impractically small. To solve this issue, a stabilization procedure was proposed in order to anticipate the effect of delays. The admittance of the emulated shunt was modified in order to obtain a closed-loop system as close as possible to an actual analog shunt.
        
        The developments were experimentally validated on a piezoelectric beam controlled by a DVA. Namely, it was shown that the expression for the maximum sampling period for stability is accurate, and that the stabilization procedure leads to a controlled system which behaves similarly to a piezoelectric structure with a shunt.
        
        A possible extension to this work could be to use the $z$-transform in order to analyze more rigorously the dynamics of the sampled-data system for large sampling periods.
    
   \section*{Funding} 
This work was supported by the SPW [WALInnov grant 1610122].

\bibliographystyle{IEEEtran}
\bibliography{biblio.bib}

\end{document}